\documentclass[12pt]{amsart}
\newtheorem{thm}{Theorem}[section]

\newtheorem{lem}[thm]{Lemma}
\newtheorem{prop}[thm]{Proposition}
\theoremstyle{definition}

\theoremstyle{remark}
\newtheorem{rem}[thm]{Remark}
\numberwithin{equation}{section}
\newcommand{\Real}{\mathbb R}
\newcommand{\Sphere}{\mathbb S}
\newcommand{\eps}{\varepsilon}

\begin{document}

\title[]{On the absolutely continuous spectrum of the Laplace-Beltrami
operator acting on $p$-forms
for a class of warped product metrics}%
\author{Francesca Antoci}
\address{Dipartimento di Matematica, Politecnico di Torino,
C.so Duca degli Abruzzi
24, 10129, Torino}
\email{antoci@calvino.polito.it}

\thanks{}%

%\date{^^^^^^^^^^^^^^^^^^^^^^^^^^^^^^^^^^^^^^^^^^^^^^^^^^^^^^^^^^^^^^^^^^^^^^^^$
%\dedicatory{^^^^^^^^^^^^^^^^^^^^^^^^^^^^^^^^^^^^^^^^^^^^^^^^^^^^^^^^^^^^^^^^^^$
%\commby{^^^^^^^^^^^^^^^^^^^^^^^^^^^^^^^^^^^^^^^^^^^^^^^^^^^^^^^^^^^^^^^^^^^^^^$
% ----------------------------------------------------------------
\begin{abstract} We explicitely compute the absolutely continuous
\\ spectrum
 of the
Laplace-Beltrami operator for $p$-forms for the class of warped
product metrics $d\sigma^2= y^{2a}dy^2 +
y^{2b}d\theta_{\Sphere^{N-1}}^2$, where $y$ is a boundary defining
function on the unit ball $B(0,1)$ in $\Real^N$.
\end{abstract}
\maketitle
% ----------------------------------------------------------------
\section{Introduction}

In the present paper we continue the investigation of the spectrum
of the Laplace-Beltrami operator acting on $p$-forms for a class
of warped product metrics started in \cite{Antoci2}. The
Riemannian manifolds considered in that paper were constructed as
follows: let $\overline M$ be a compact $N$-dimensional manifold
with boundary, and let $y$  be any boundary-defining function. We
endowed the interior $M$ of $\overline M$ with a Riemannian metric
$ds^2$ such that in a small tubular neighbourhood of $\partial M$
in $M$ $ds^2$ takes the form $$ds^2=e^{-2(a+1)t}dt^2 +
e^{-2bt}d\theta^2_{\partial M},$$ where $t=- \log y \in
(c,+\infty)$ and $d\theta^2_{\partial M}$ is a Riemannian metric
on $\partial M$. For $a\leq -1$, the manifold $M$ is complete,
hence the Laplace-Beltrami operator $\Delta^p_M$ is essentially
selfadjoint on the smooth compactly supported $p$-forms. In
\cite{Antoci2}, we computed explicitely, with the exception of the
point $0$, the essential spectrum of $\Delta^p_M$ in dependence on
the parameters $a$ and $b$. Moreover, under the assumption of
rotational symmetry, that is, assuming that $\partial M=
\Sphere^{N-1}$, we were able to check the belonging of $0$ to the
essential spectrum of $\Delta^p_M$ , and hence to achieve a
complete description of the essential spectrum.
\par In the present paper, instead, we are concerned with the
absolutely continuous (and, partly, with the singularly
continuous) spectrum.\par  In \cite{Eichhorn}, Eichhorn  showed
that the essential spectrum of $\Delta^p_M$ coincides with the
essential spectrum of the Friedrichs extension $(\Delta^p_M)^F$ of
the restriction of $\Delta^p_M$ to the smooth $p$-forms with
compact support in $(c,+\infty)\times
\partial M$. Hence, in order to achieve the results in
\cite{Antoci2}, it sufficed to know the behaviour of the
Riemannian metric only in a tubular neighbourhood of the boundary.
\par
As for the absolutely continuous (and the singularly continuous)
spectrum of $\Delta^p_M$, instead, to our knowledge no result of
the sort of \cite{Eichhorn} is available. As a consequence, in
order to compute these parts of the spectrum, we need global
information on the Riemannian manifold.
\par In the present paper we restrict our
attention to the case in which $M$ is the unitary open ball
$B(0,1)$ in $\Real^N$ endowed with a Riemannian metric $ds^2$
given by
\begin{equation}\label{metric0}ds^2:=f(t)dt^2 +g(t)d\theta_{\Sphere^{N-1}}^2,
\end{equation}
where $t=2 {\rm settanh}(\| x\|)$, $f(t)>0$, $g(t)>0$ for every
$t\in (0,+\infty)$, and $d\theta_{\Sphere^{N-1}}^2$ is the
standard Riemannian metric on $\Sphere^{N-1}$. Moreover, we assume
that $f(t)=1$ and $g(t)=t^2$ for $0<t<\epsilon$, whilst
$f(t)=e^{-2(a+1)t}$ and $g(t)=e^{-2bt}$ for $t>c>\epsilon$.
\par On one hand, these assumptions give us a complete knowledge
of the essential spectrum (see \cite{Antoci2}); on the other hand,
they let us employ the radial decomposition techniques developed
by Dodziuk (\cite{Dodziuk}), Donnelly (\cite{Donnelly}) and
Eichhorn (\cite{Eichhorn}, \cite{Eichhorn2}). The decomposition
consists of two steps: first, thanks to the Hodge decomposition on
$\Sphere^{N-1}$, we write any $p$-form $\omega$ on $M$ as
\begin{equation}\label{dec}\omega= \omega_{1\delta}\oplus
\omega_{2d}\wedge dt \oplus (\omega_{1d}\oplus
\omega_{2\delta}\wedge dt),\end{equation} where $\omega_{1\delta}$
(resp. $\omega_{1d}$) is a coclosed (resp. closed) $p$-form on
$\Sphere^{N-1}$ parametrized by $t$ and $\omega_{2\delta}$ (resp.
$\omega_{2d}$) is a coclosed (resp. closed) $(p-1)$-form on
$\Sphere^{N-1}$ parametrized by $t$. This gives the orthogonal
decomposition $$L^2_p(M)= {\mathcal L}_{p,1}(M) \oplus {\mathcal
L}_{p,2}(M)\oplus {\mathcal L}_{p,3}(M),$$ and, since $\Delta^p_M$
is invariant, the corresponding decomposition $$\Delta^p_M=
\Delta^p_{M1}\oplus \Delta^p_{M2}\oplus \Delta^p_{M3}.$$ Since
 $$
\sigma_{\rm ac}(\Delta^p_M)= \bigcup_{i=1}^3\sigma_{\rm
ac}(\Delta^p_{Mi}),$$$$ \sigma_{\rm sc}(\Delta^p_M)=
\bigcup_{i=1}^3\sigma_{\rm sc}(\Delta^p_{Mi}),$$ we can reduce
ourselves to the study of the absolutely continuous (and of the
singularly continuous) spectrum of $\Delta^p_{Mi}$ for $i=1,2,3$.
\par The second step consists of the decomposition of
$\omega_{1\delta}$ (resp. of $\omega_{2d}$, $\omega_{2\delta}$)
according to an orthonormal basis of coclosed $p$-eigenforms
(resp. closed and coclosed $(p-1)$-eigenforms) of
$\Delta^p_{\Sphere^{N-1}}$ (resp. of
$\Delta^{p-1}_{\Sphere^{N-1}}$) on $\Sphere^{N-1}$. In this way,
up to a unitary equivalence, the spectral analysis of
$\Delta^p_{Mi}$, for $i=1,2,3$, can be reduced to the
investigation of the spectra of a countable number of selfadjoint
Sturm-Liouville operators $D_{i\lambda_k}$ on the half-line
$(0,+\infty)$, parametrized by the eigenvalues $\lambda_k^p$,
$k\in \mathbb N$ of $\Delta^p_{\Sphere^{N-1}}$ on $\Sphere^{N-1}$
if $i=1$, and by the eigenvalues $\lambda_k^{p-1}$ of
$\Delta^{p-1}_{\Sphere^{N-1}}$ if $i=2,3$. In particular, we have
that for $i=1,2,3$ $$ \sigma_{\rm ac}(\Delta^p_{Mi})=
\bigcup_{k\in \mathbb N}\sigma_{\rm ac}(D_{i\lambda_k}),$$$$
\sigma_{\rm sc}(\Delta^p_{Mi})= \bigcup_{k\in \mathbb
N}\sigma_{\rm sc}(D_{i\lambda_k}).$$
 \par Actually, since the Hodge
$*$ operator maps isometrically $p$-forms of ${\mathcal
L}_{p,1}(M)$ into $(N-p)$-forms of ${\mathcal L}_{N-p,2}(M)$ and
viceversa, it suffices to consider the cases $i=1,3$. Moreover, it
turns out that, since the absolutely continuous spectrum of
$\Delta^p_M$ is contained in the essential spectrum of
$\Delta^p_M$, which we know from \cite{Antoci2}, in order to
compute the absolutely continuous spectrum of $\Delta^p_M$ it
suffices to study the absolutely continuous spectrum of
$D_{1\lambda_k^p}$ for any $k\in \mathbb N$: indeed, for any
$a\leq -1$, $b\in \Real$, $p\in [0,N]$, we find that
$\bigcup_{k\in \mathbb N}\sigma_{\rm
ac}(D_{1\lambda_k^p})=\sigma_{\rm ess}(\Delta^p_M).$
\par The
absolutely continuous spectrum (and the singularly continuous
spectrum) of the operators $D_{1\lambda_k^p}$ is computed through
perturbation theory; on one hand, this required a subtle
investigation of their domains. On the other hand, since the
operators $D_{1\lambda_k^p}$ act on the one-dimensional half-line
$(0,+\infty)$ and have strongly divergent potential terms at zero,
in order to study their spectra we had to prove modified versions
of the classical Agmon-Kato-Kuroda Theorem (\cite{Reed-Simon4})
and Lavine Theorem (\cite{Lavine}). In particular, we had to
choose properly the unperturbed operators employed in the
perturbation techniques.
\par Let us briefly discuss our results. If $a=-1$, $b<0$, the
situation is similar to the hyperbolic case; we find that for
$0\leq p \leq N$ $$\sigma_{\rm ac}(\Delta^p_M)=\left[ \min
\left\{\left(\frac{N-2p-1}{2}\right)^2b^2,
\left(\frac{N-2p+1}{2}\right)^2b^2\right\},+\infty\right).$$ If
$a\leq -1$ and $b=0$, for $0\leq p \leq N$ $$\sigma_{\rm
ac}(\Delta^p_M)= [\overline \lambda_p, +\infty),$$ where
$\overline \lambda_p$ is the minimum between the lowest eigenvalue
$\lambda_0^p$ of $\Delta^p_{\Sphere^{N-1}}$ and the lowest
eigenvalue $\lambda_0^{p-1}$ of
$\Delta^{p-1}_{\Sphere^{N-1}}$.\par For $a=-1$ and $b>0$, if
$1<p<N-1$ $\sigma_{\rm ac}(\Delta^p_M)=\emptyset$, whilst if $p\in
\left\{ 0,1,N-1,N\right\}$ $\sigma_{\rm ac}(\Delta^p_M)=
\left[\left(\frac{N-1}{2}\right)^2b^2,+\infty\right)$. \par If
$a<-1$ and $b<0$, for $0\leq p \leq N$ $\sigma_{\rm
ac}(\Delta^p_M)=[0,+\infty)$; finally, if $a<-1$ and $b>0$, for
$1<p<N-1$ $\sigma_{\rm ac}(\Delta^p_M)= \emptyset$, whilst for
$p\in \left\{ 0,1,N-1,N \right\}$ $\sigma_{\rm
ac}(\Delta^p_M)=[0,+\infty)$.\par As for the singularly continuous
spectrum, in any case we found that $\sigma_{\rm sc}(\Delta^p_M)=
\sigma_{\rm sc}(\Delta^p_{M3})$, whilst $\sigma_{\rm
sc}(\Delta^p_{M1})= \sigma_{\rm sc}(\Delta^p_{M2})= \emptyset$.
\par
It would be interesting to complete the analysis of the spectrum
of $\Delta^p_M$, computing its singularly continuous spectrum.
This problem can be reduced to the determination of the singularly
continuous spectrum of $D_{3\lambda_k^{p-1}}$ for any $k\in
\mathbb N$; this turns out to be a hard task because
$D_{3\lambda_k^{p-1}}$ is a coupled system of Sturm-Liouville
operators on the half-line $(0,+\infty)$ with strongly divergent
potentials at zero, for which the application of perturbation
techniques is difficult.\par The paper is organized as follows: in
section 2 we introduce some preliminary facts and basic notations.
In section 3 we describe in some detail the decomposition
techniques. The calculus of the absolutely continuous spectrum
(and, partly, of the singularly continuous spectrum) of
$\Delta^p_M$ is performed in section 4 for $a=-1$ and in section 5
for $a<-1$.

\section{Preliminary facts and notations}
For $N\geq 2$, let $\overline{B}(0,1)$ denote the closed unit ball
$$\overline{B}(0,1)=\left\{  x= (x_1,...,x_N) \in \Real^N\,|\,
x_1^2+...+x_N^2\leq 1\right\}, $$ and let $\Sphere^{N-1}$ denote
the sphere $$\Sphere^{N-1}= \left\{(x_1,...,x_N)\in \Real^N\,|\,
x_1^2+...+x_N^2=1 \right\} ,$$ endowed with a coordinate system
$(U_i, \Theta_i)$, $i=2,...,k+1$, $\Theta_i: U_i \rightarrow
\Real^{N-1}$.\par Let us consider the interior of
$\overline{B}(0,1)$, $$B(0,1)=\left\{(x_1,...,x_N)\in \Real^N\,|\,
x_1^2+...+x_N^2<1 \right\}, $$ with the coordinate system
$(V_i,\Phi_i)$, for $i=1,...,k+1$, defined in the following way:
in a neighbourhood of $0$, for some $\delta >0$, $$ V_1= \left\{
(x_1,...,x_N)\in \Real^N\,|\, x_1^2 +...+x_N^2 < \delta\right\}$$
and $$ \Phi_1(x_1,...,x_N)=(x_1,...,x_N),$$ whilst for $i>1$, $ x
\not=0$, $$V_i=\left\{  x \in \Real^N \,|\, \frac{ x}{\| x\|}\in
U_i\right\},$$ $$ \Phi_i: V_i \longrightarrow (0,+\infty) \times
\Theta_i(U_i),$$ $$\Phi_i(x_1,...,x_N)= \left(2\, {\rm settanh}
(\| x\|), \Theta_i\left(\frac{ x}{\|
x\|}\right)\right)=:(t,\theta_i).$$ We denote by $M$ the manifold
$B(0,1)$, endowed with a Riemannian metric $ds^2$ such that on
$\Phi_i(V_i)$, for $i>1$,
\begin{equation}\label{metric} ds^2:=f(t) dt^2 + g(t)
d\theta_{\Sphere^{N-1}}^2,\end{equation} where $f(t)>0$, $g(t)>0$
for every $t\in (0,+\infty)$ and $d \theta_{\Sphere^{N-1}}^2$ is
the standard metric on $\Sphere^{N-1}$. $ds^2$ is well-defined on
$B(0,1) \setminus \left\{ 0\right\}$.  \par We suppose that for $t
> c>0$, $a\in \Real$, $b\in \Real$
\begin{equation}\label{condinf} f(t)= e^{-2(a+1)t}, \quad \quad
g(t)=e^{-2bt}.\end{equation} As for the behaviour as
$t\rightarrow 0$, we suppose that for $t\in (0,\epsilon)$
($\epsilon= 2\,{\rm settanh} (\delta)$)
\begin{equation}\label{condzero} f(t)\equiv 1, \quad \quad g(t)=t^2.
\end{equation}
This assures that $ds^2$ can be extended to a smooth Riemannian
metric on the whole manifold $M$; indeed, for $t\in (0,\epsilon)$,
$ds^2$ is the expression, in polar coordinates, of the Euclidean
metric on $\Real^N$. \par It is well-known (see \cite{Melrose})
that a Riemannian metric of this kind is complete if and only if
$a\leq -1$. Hence throughout the paper we will suppose that $a\leq
-1$.
\par For $p=0,...,N$, we will denote by
$C^{\infty}(\Lambda^p(M))$ the space of all smooth $p$-forms on
$M$, and by $C^{\infty}_c(\Lambda^p(M))$ the set of all smooth,
compactly supported $p$-forms on $M$. For any $\omega\in
C^{\infty}(\Lambda^p(M))$, we will denote by
$|\omega(t,\theta)|_M$ the norm induced by the Riemannian metric
on the fiber over $(t,\theta)$, given in local coordinates by
$$|\omega(t,\theta)|_M^2=g^{i_1j_1}(t,\theta)...g^{i_pj_p}(t,\theta)
\omega_{i_1...i_p}(t,\theta)\omega_{j_1...j_p}(t, \theta),$$ where
$g^{ij}$ is the expression of the Riemannian metric in local
coordinates. We will denote by $d^p_M$, $*_M$, $\delta^p_M$,
respectively, the differential, the Hodge $*$ operator and the
codifferential on $M$, defined as in \cite{deRham}. $\Delta^p_M$
will stand for the Laplace-Beltrami operator acting on $p$-forms
$$ \Delta^p_M=d^{p-1}_M\delta^p_M+\delta^{p+1}_M d^p_M,$$ which is
expressed in local coordinates by the Weitzenb\"ock formula
$$((\Delta^p_M)\omega)_{i_1...i_p}= - g^{ij}\nabla_i \nabla_j
\omega_{i_1...i_p}+\sum_j R^{\alpha}_{j}
\omega_{i_1...\alpha...i_p} + \sum_{j,l\not=j}
R^{\alpha\,\beta\,}_{\,i_j\,i_l} \omega_{\alpha
i_1...\beta...i_p}, $$ where $\nabla_i\omega$ is the covariant
derivative of $\omega$ with respect to the Riemannian metric, and
$R^i_j$, $R^{i\,j\,}_{\,k\,l}$ denote respectively the local
components of the Ricci tensor and the Riemann tensor induced by
the Riemannian metric. As usual, $L^2_p(M)$ will denote the
completion of $C^{\infty}_c(\Lambda^p(M))$ with respect to the
norm $\|\omega\|_{L^2_p(M)}$ induced by the scalar product
$$\langle\omega,\tilde\omega\rangle_{L^2_p(M)}:= \int_M \omega
\wedge
*_M\tilde \omega
; $$
$\|\omega\|_{L^2_p(M)}$ reads also $$ \|\omega\|^2_{L^2_p(M)}=
\int_M |\omega(t,\theta)|_M^2 dV_M,$$ where $dV_M$ is the volume
element of $(M,ds^2)$.\par It is well-known that, since the
Riemannian metric on $M$ is complete, the Laplace-Beltrami
operator is essentially selfadjoint on
$C^{\infty}_c(\Lambda^p(M))$, for $p=0,...,N$. We will denote by
$\Delta^p_M$ also its closure.\par Let us end this section with
some notations and preliminary facts in spectral analysis. If $H$
is any selfadjoint operator acting in a Hilbert space ${\mathcal
H}$, $$ H:{\mathcal D}(H)\subseteq {\mathcal H} \longrightarrow
{\mathcal H},$$ we will denote by $\sigma_{\rm ess}(H)$ the
essential spectrum of $H$, that is, the spectrum of $H$ minus the
isolated eigenvalues of finite multiplicity. Following
\cite{Kato}, $E_H(\mu)$ ($\mu \in \Real$) will stand for the
spectral family associated to the operator $H$; moreover,
$P_H(\mu)$ will denote the projection $E_H(\mu)\ominus E_H(\mu-0)$
(where $E_H(\mu -0)= {\rm s}-\lim_{\epsilon \rightarrow
0}E_H(\mu)$), whilst $E_H(S)$ will stand for the spectral measure
of any Borel set $S \subseteq \Real$. As usual, we will denote by
${\mathcal H}_p(H)$ the closed subset of ${\mathcal H}$ spanned by
all the eigenfunctions of $H$, and by ${\mathcal H}_c(H)$ its
orthogonal complement in ${\mathcal H}$; correspondingly, we will
denote by $\sigma_p(H)$ the set of all the eigenvalues of $H$ and
by $\sigma_c(H)$ the spectrum of the restriction of $H$ to
${\mathcal H}_c(H)$. Following \cite{Reed-Simon1}, we will denote
by ${\mathcal H}_{ac}(H)$ the subset of absolute continuity of
$H$, defined as the set of all $u \in {\mathcal H}$ such that
$\langle E_H(S)u,u \rangle_{\mathcal H}=0$ for any Borel set $S$
whose  Lebesgue measure $|S|$ is equal to zero. ${\mathcal
H}_{sc}(H)$ will stand for the set ${\mathcal H}_c(H)\ominus
{\mathcal H}_{ac}(H)$. Accordingly, we will denote by $\sigma_{\rm
ac}(H)$ (resp. $\sigma_{\rm sc}(H)$) the absolutely (resp.
singularly) continuous spectrum of $H$, defined as the spectrum of
the restriction of $H$ to the subspace ${\mathcal H}_{ac}(H)$
(resp. ${\mathcal H}_{sc}(H)$). \par Finally, let us recall the
following basic facts, whose proof is elementary and is therefore
omitted:
\begin{lem}\label{basic1} Let $H$ be a selfadjoint operator acting on a
Hilbert space ${\mathcal H}$, $H:D(H)\subseteq {\mathcal H} \rightarrow
{\mathcal H}$. Then
\begin{enumerate}
\item if $\mu \in \Real$ is an isolated eigenvalue of $H$, then $\mu
\notin \sigma_{\rm ac}(H)$ (resp. $\mu \notin \sigma_{\rm sc}(H)$); as a
consequence, $\sigma_{\rm ac}(H)\subseteq \sigma_{\rm ess}(H)$ (resp.
$\sigma_{\rm sc}(H)\subseteq \sigma_{\rm ess}(H)$); \\
\item if ${\mathcal H}= \oplus_{k \in \mathbb N} {\mathcal H}_k$, where
${\mathcal H}_k$, for every $k\in \mathbb N$, is a closed
subspace of ${\mathcal H}$
(possibly empty), and if $H$ splits accordingly as $H= \oplus_{k\in
\mathbb
N}H_k$, where for every $k\in \mathbb N$ $H_k= H_{|{\mathcal H}_k}$, then
$\sigma_{\rm ac}(H)= \bigcup_{k\in \mathbb N} \sigma_{\rm ac}(H_k)$ and
$\sigma_{\rm sc}(H)= \bigcup_{k\in \mathbb N} \sigma_{\rm sc}(H_k)$;
 \\
\item for any constant $K\in \Real$, ${\mathcal
H}_{ac}(H+K)={\mathcal H}_{ac}(H)$
(resp. ${\mathcal H}_{sc}(H+K)={\mathcal H}_{sc}(H)$); as
a consequence, $\sigma_{\rm ac}(H+K)=\sigma_{\rm ac}(H)+K$ (resp.
$\sigma_{\rm sc}(H+K)=\sigma_{\rm sc}(H)+K$). \end{enumerate}
\end{lem}

\section{Hodge decomposition}
In the present section let us suppose that the Riemannian metric
$ds^2$ in $(0,+\infty)\times \Sphere^{N-1}$ takes the form
\begin{equation}\label{fg}ds^2= f(t)\,dt^2 +
g(t)\,d\theta_{\Sphere^{N-1}}^2,\end{equation} where $f(t)>0$ and
$g(t)>0$ for any $t\in (0,+\infty)$. \par Given $\omega \in
C^{\infty}(\Lambda^p(M))$, let us write, for $(t,\theta)\in
(0,+\infty)\times \Sphere^{N-1}$
\begin{equation}\label{decom1} \omega(t,\theta)= \omega_1(\theta) +
\omega_2(\theta)
\wedge dt,\end{equation}
where $\omega_1$ and $\omega_2$ are respectively a $p$-form and a
$(p-1)$-form on $\Sphere^{N-1}$ depending on $t$. An easy
computation shows that $*_M \omega$ can be expressed in terms of
decomposition
(\ref{decom1}) as
\begin{multline}\label{Hodge1}*_M \omega=
(-1)^{N-p}g^{\frac{N-2p+1}{2}}(t) f^{-\frac{1}{2}}(t)
*_{\Sphere^{N-1}} \omega_2 \\+
g^{\frac{N-2p-1}{2}}(t)f^{\frac{1}{2}}(t)*_{\Sphere^{N-1}}
\omega_1\wedge dt, \end{multline} where $*_{\Sphere^{N-1}}$
denotes the Hodge $*$ operator on $\Sphere^{N-1}$. Moreover,
$d^p_M$ and $\delta^p_M$ split respectively as
\begin{equation}\label{dM1}d^p_M\omega= d^p_{\Sphere^{N-1}}\omega_1 + \left\{
(-1)^p \frac{\partial \omega_1}{\partial t}+
d^{p-1}_{\Sphere^{N-1}} \omega_2 \right\}\wedge dt ,\end{equation}
\begin{multline}\label{deltaM1}\delta^p_M \omega= g^{-1}(t)
\delta^p_{\Sphere^{N-1}} \omega_1 + (-1)^p f^{-\frac{1}{2}}
g^{\frac{-N-1+2p}{2}}\frac{\partial}{\partial t}
\left(f^{-\frac{1}{2}} g^{\frac{N+1-2p}{2}}\omega_2\right)\\+
g^{-1} \delta^{p-1}_{\Sphere^{N-1}} \omega_2 \wedge
dt,\end{multline} where $p$ is the degree of $\omega$,
$d^p_{\Sphere^{N-1}}$ is the differential on $\Sphere^{N-1}$ and
$\delta^{p-1}_{\Sphere^{N-1}}$ is the codifferential on
$\Sphere^{N-1}$.
\par Moreover, the $L^2$-norm of $\omega \in
C^{\infty}(\Lambda^p(M))\cap L^2_p(M)$ can be written as
\begin{multline}\label{norm1} \|\omega\|^2_{L^2_p(M)}= \int_0^{+\infty}
g^{\frac{N-2p-1}{2}}(s) f^{\frac{1}{2}}(s)
\|\omega_1(s)\|^2_{L^2_p(\Sphere^{N-1})} \,ds \\ +
\int_0^{+\infty} g^{\frac{N+1-2p}{2}}(s)f^{-\frac{1}{2}}(s)
\|\omega_2(s)\|^2_{L^2_{p-1}(\Sphere^{N-1})} \,ds,
\end{multline}
where $\|.\|_{L^2_p(\Sphere^{N-1})}$ is the $L^2$-norm for
$p$-forms on $\Sphere^{N-1}$.   \par From (\ref{dM1}) and
(\ref{deltaM1}), a lengthy but straightforward computation gives
$$ \Delta^p_M \omega= (\Delta^p_M \omega)_1 + (\Delta^p_M
\omega)_2\wedge dt,$$ where
\begin{multline}\label{LB11}(\Delta^p_M \omega)_1=
g^{-1}(t)\Delta^p_{\Sphere^{N-1}} \omega_1 + (-1)^p f^{-1}(t)
g^{-1}(t) \frac{\partial g}{\partial t}
d^{p-1}_{\Sphere^{N-1}}\omega_2
\\ -f^{-\frac{1}{2}}(t) g^{\frac{-N+1+2p}{2}}(t)
\frac{\partial}{\partial t} \left( f^{-\frac{1}{2}}(t)
g^{\frac{N-1-2p}{2}}(t) \frac{\partial \omega_1}{\partial t}
\right)
\end{multline} and
\begin{multline}\label{LB12}(\Delta^p_M \omega)_2 =
g^{-1}(t) \Delta^{p-1}_{\Sphere^{N-1}} \omega_2  +
(-1)^{p}g^{-2}(t) \frac{\partial g}{\partial t}
\delta^p_{\Sphere^{N-1}} \omega_1  \\ - \frac{\partial}{\partial
t} \left\{f^{-\frac{1}{2}}(t) g^{\frac{-N-1+2p}{2}}(t)
\frac{\partial}{\partial t}\left(f^{-\frac{1}{2}}(t)
g^{\frac{N+1-2p}{2}}(t) \omega_2\right) \right\}. \end{multline}
Here we denote by $\Delta^p_{\Sphere^{N-1}}$ the Laplace-Beltrami
operator on $\Sphere^{N-1}$.\par Since for every smooth $\omega\in
L^2_p(M)$ we have that $\omega_1 \in L^2_p(M)$, $\omega_2 \wedge
dt \in L^2_p(M)$ and $$ \langle\omega_1,\omega_2\wedge
dt\rangle_{L^2_p(M)}=0,$$ (\ref{decom1}) gives rise to an
orthogonal decomposition of $L^2_p(M)$ into two closed subspaces.
However, (\ref{LB11}) and (\ref{LB12}) show that $\Delta^p_M$ is
not invariant under this decomposition, and further decompositions
are required.\par It is well-known that, for $0\leq p \leq N-1$,
$$ C^{\infty}(\Lambda^p(\Sphere^{N-1}))=
dC^{\infty}(\Lambda^{p-1}(\Sphere^{N-1}))\oplus \delta
C^{\infty}(\Lambda^{p+1}(\Sphere^{N-1})) \oplus {\mathcal
H}^p(\Sphere^{N-1}),$$ where ${\mathcal H}^p(\Sphere^{N-1})$ is
the space of harmonic $p$-forms on $\Sphere^{N-1}$, and the
decomposition is orthogonal in $L^2_p(\Sphere^{N-1})$. Hence, for
$1\leq p \leq N-1$, every $\omega\in L^2_p(M)\cap
C^{\infty}(\Lambda^p(M))$ can be written as $$\omega=
\omega_{1\delta}\oplus \omega_{2d}\wedge dt \oplus (\omega_{1d}
\oplus \omega_{2\delta}\wedge dt),$$ where $\omega_{1\delta}$
(resp. $\omega_{1d}$) is a coclosed (resp. closed) $p$-form on
$\Sphere^{N-1}$ parametrized by $t$, and $\omega_{2\delta}$ (resp.
$\omega_{2d}$) is a coclosed (resp. closed) $(p-1)$-form on
$\Sphere^{N-1}$ parametrized by $t$. By closure, we get the
orthogonal decomposition $$ L^2_p(M)= {\mathcal L}_{p,1}(M)\oplus
{\mathcal L}_{p,2}(M) \oplus{\mathcal L}_{p,3}(M),$$ where, for
every $\omega \in L^2_p(M)\cap C^{\infty}(\Lambda^p(M))$,
$$\omega_{1\delta}\in {\mathcal L}_{p,1}(M),$$ $$
\omega_{2d}\wedge dt \in {\mathcal L}_{p,2}(M)$$ and $$
(\omega_{1d} \oplus \omega_{2\delta}\wedge dt)\in {\mathcal
L}_{p,3}(M).$$ Since $$d^p_{\Sphere^{N-1}}
\Delta^p_{\Sphere^{N-1}} =\Delta^{p+1}_{\Sphere^{N-1}}
d^p_{\Sphere^{N-1}},\quad \quad \delta^p_{\Sphere^{N-1}}
\Delta^p_{\Sphere^{N-1}} =\Delta^{p-1}_{\Sphere^{N-1}}
\delta^p_{\Sphere^{N-1}},$$ $$\frac{\partial}{\partial
t}d^p_{\Sphere^{N-1}}=d^p_{\Sphere^{N-1}} \frac{\partial}{\partial
t}, \quad \quad \frac{\partial}{\partial
t}\delta^p_{\Sphere^{N-1}}=\delta^p_{\Sphere^{N-1}}
\frac{\partial}{\partial t},$$ the Laplace-Beltrami operator is
invariant under this decomposition, and $$\Delta^p_M=\Delta^p_{M1}
\oplus \Delta^p_{M2} \oplus \Delta^p_{M3},$$ where, for $i=1,2,3$,
$\Delta^p_{Mi}$ is the restriction of $\Delta^p_M$ to ${\mathcal
L}_{p,i}(M)$. We remark that, for $i=1,2,3$, $\Delta^p_{Mi}$ is
essentially selfadjoint on $C^{\infty}_c(\Lambda^p(M))\cap
{\mathcal L}_{p,i}(M)$. In view of Lemma \ref{basic1}, for $1\leq
p \leq N-1$, $$ \sigma_{\rm ac}(\Delta^p_M)= \bigcup_{i=1}^3
\sigma_{\rm ac}(\Delta^p_{Mi}),$$ $$ \sigma_{\rm
sc}(\Delta^p_M)=\bigcup_{i=1}^3 \sigma_{\rm sc} (\Delta^p_{Mi}).$$
For $p=0$ (resp. $p=N$), any $\omega \in L^2_p(M)$ can be written
as $\omega=\omega_{1\delta}$ (resp. $\omega= \omega_{2d} \wedge
dt$), where $\omega_{1\delta}$ (resp. $\omega_{2d}$) is a coclosed
(resp. closed) $0$-form (resp. $(N-1)$-form) parametrized by $t$
on $\Sphere^{N-1}$. Hence $L^2_0(M)= {\mathcal L}_{0,1}(M)$ (resp.
$L^2_N(M)= {\mathcal L}_{N,2}(M)$), and $\Delta^0_M=
\Delta^0_{M1}$ (resp. $\Delta^N_M=\Delta^N_{M2}$). Again,
$\Delta^0_{M1}$ (resp. $\Delta^N_{M2}$) is essentially selfadjoint
on $C^{\infty}_c(\Lambda^0(M))\cap {\mathcal L}_{0,1}(M)$ (resp.
on $C^{\infty}_c(\Lambda^N(M))\cap {\mathcal L}_{N,2}(M)$).
\par Hence, for every $p\in [0,N]$, in order to
determine the spectral properties of $\Delta^p_M$ it suffices to
study the corresponding properties of $\Delta^p_{Mi}$, $i=1,2,3$.
\par To this purpose, let us introduce further decompositions.
First of all, for any $\omega \in L^2_p(M)\cap
C^{\infty}(\Lambda^p(M))$ we decompose $\omega_{1\delta}$
according to an orthonormal basis $\left\{\tau_{1k}\right\}_{k \in
\mathbb N}$ of coclosed eigenforms of $\Delta^p_{\Sphere^{N-1}}$;
this yields
\begin{equation}\label{1delta} \omega_{1\delta}=\oplus_k
h_k(t) \tau_{1k},\end{equation} where $h_k(t)\tau_{1k}\in
L^2_p(M)$ for every $k \in \mathbb N$, and the sum is orthogonal
in $L^2_p(M)$, thanks to (\ref{metric}). By closure, we get the
decomposition $$ {\mathcal L}_{p,1}(M)= \oplus_{k \in \mathbb
N}{\mathcal L}_{p,1,k}(M),$$ where for $\omega \in L^2_p(M) \cap
C^{\infty}(\Lambda^p(M))$ $$h_k(t)\tau_{1k}\in {\mathcal
L}_{p,1,k}(M)$$ for every $k\in \mathbb N$. We will call $p$-form
of type I any $p$-form $\omega\in L^2_p(M)$ such that $\omega \in
{\mathcal L}_{p,1,k}(M)$ for some $k\in \mathbb N$.
\par
For every $k \in \mathbb N$, let us denote by $\lambda_k^p$ the
eigenvalue associated to $\tau_{1k}$. Since for every $k\in
\mathbb N$
\begin{multline}\label{h1}\Delta^p_{M1}(h(t)\tau_{1k})=
\left\{ \frac{\lambda^p_k}{g(t)} \right. \\ \left. -
f(t)^{-\frac{1}{2}}g(t)^{\frac{-N+1+2p}{2}}\frac{\partial}{\partial
t}\left(f(t)^{-\frac{1}{2}}g(t)^{\frac{N-1-2p}{2}}\frac{\partial
h}{\partial t}\right)\right\} \tau_{1k}, \end{multline}
$\Delta^p_{M1}$ is invariant under decomposition (\ref{1delta}).
Moreover, if $\omega=h(t)\tau_{1k}$ $$\|\omega\|^2_{L^2_p(M)}=
\int_0^{+\infty}g(s)^{\frac{N-2p-1}{2}}f(s)^{\frac{1}{2}}h(s)^2\,ds.$$
Thus, $\Delta^p_{M1}$ is unitarily equivalent to the direct sum
over $k \in \mathbb N$ of certain selfadjoint operators
$\Delta_{1\lambda^p_k}$ in $L^2((0,+\infty),
g^{\frac{N-2p-1}{2}}f^{\frac{1}{2}})$ such that  $$C^{\infty}_c
(0,+\infty) \subseteq {\mathcal D}(\Delta_{1\lambda^p_k})$$ and
for every $h \in C^{\infty}_c (0,+\infty)$
\begin{multline} \Delta_{1\lambda^p_k}h=
\frac{\lambda^p_k}{g(t)} h(t)  -
f(t)^{-\frac{1}{2}}g(t)^{\frac{-N+1+2p}{2}}\frac{\partial}{\partial
t}\left(f(t)^{-\frac{1}{2}}g(t)^{\frac{N-1-2p}{2}}\frac{\partial
h}{\partial t}\right).\end{multline} If we set
\begin{equation}\label{trasf1}w(t)=
h(t)f(t)^{\frac{1}{4}}g(t)^{\frac{N-2p-1}{4}},
\end{equation}
a direct (but lengthy) computation shows that $\Delta^p_{M1}$ is
unitarily equivalent to the direct sum, over $k\in \mathbb N$, of
some selfadjoint operators $D_{1\lambda^p_k}$ in $L^2(0,+\infty)$
such that $$C^{\infty}_c(0,+\infty)\subseteq {\mathcal
D}(D_{1\lambda^p_k})$$ and for every $w \in
C^{\infty}_c(0,+\infty)$
\begin{multline}\label{w1} D_{1\lambda^p_k}w = -
\frac{\partial}{\partial t}\left( \frac{1}{f} \frac{\partial
w}{\partial t} \right)+ \left\{ -\frac{7}{16} \frac{1}{f^3}
\left(\frac{\partial f}{\partial t}\right)^2+ \frac{1}{4}
\frac{1}{f^2} \frac{\partial^2 f}{\partial t^2}\right.\\ -
\frac{1}{2} \frac{1}{f^2}\frac{\partial f}{\partial
t}\frac{(N-1-2p)}{4} \frac{1}{g}\frac{\partial g}{\partial t} +
\frac{1}{f} \frac{(N-2p-1)}{4}
\frac{(N-2p-5)}{4}\frac{1}{g^2}\left(\frac{\partial g}{\partial
t}\right)^2\\ \left. +
\frac{1}{f}\frac{(N-2p-1)}{4}\frac{1}{g}\frac{\partial^2 g
}{\partial t^2} + \frac{\lambda^p_k}{g}\right\} w.
\end{multline}
Analogously, for every $\omega \in L^2_p(M)\cap
C^{\infty}(\Lambda^p(M))$ we decompose $\omega_{2d}$ according to
an orthonormal basis of closed eigenforms $\left\{\tau_{2k}
\right\}_{k \in \mathbb N}$ of $\Delta^{p-1}_{\Sphere^{N-1}}$:
\begin{equation}\label{2d} \omega_{2d}\wedge dt=\oplus_k
h_k (t) \tau_{2k}\wedge dt.\end{equation} Correspondingly, by
closure we get the orthogonal decomposition $${\mathcal
L}_{p,2}(M)= \oplus_{k\in \mathbb N}{\mathcal L}_{p,2,k}(M);$$ we
will call $p$-form of type II a $p$-form $\omega\in L^2_p(M)$ such
that $\omega \in {\mathcal L}_{p,2,k}(M)$ for some $k\in \mathbb
N$.
\par For every $k \in \mathbb N$ $$\Delta^p_{M2}(h(t)\tau_{2k}\wedge
dt )= (\Delta_{2\lambda^{p-1}_k} h )\tau_{2k}\wedge dt,$$ where
\begin{multline}\label{h2}
\Delta_{2\lambda^{p-1}_k} h = \frac{\lambda^{p-1}_k}{g(t)} h(t)
\\ -\frac{\partial}{\partial t}\left\{
f(t)^{-\frac{1}{2}}g(t)^{\frac{-N-1+2p}{2}}
\frac{\partial}{\partial t}\left( f(t)^{-\frac{1}{2}}
g(t)^{\frac{N+1-2p}{2}}h(t) \right) \right\}. \end{multline} Here,
again, for every $k \in \mathbb N$ we denote by $\lambda_k^{p-1}$
the eigenvalue of $\Delta^{p-1}_{\Sphere^{N-1}}$ corresponding to
the eigenform $\tau_{2k}$. \par If $\omega= h(t) \tau_{2k}\wedge
dt$, then $$\|\omega\|^2_{L^2_p(M)}=
\int_0^{+\infty}g(s)^{\frac{N-2p+1}{2}}f(s)^{-\frac{1}{2}}h(s)^2\,ds.$$
Thus, if we set
\begin{equation}\label{trasf2} w(t)= h(t)
f(t)^{-\frac{1}{4}}g(t)^{\frac{N+1-2p}{4}},
\end{equation}
we find that $\Delta^p_{M2}$ is unitarily equivalent to the direct
sum, over $k \in \mathbb N$, of certain selfadjoint operators
$D_{2 \lambda^{p-1}_k}$ in $L^2(0,+\infty)$ such that
$$C^{\infty}_c(0,+\infty) \subseteq {\mathcal D} (D_{2
\lambda^{p-1}_k})$$ and for every $w\in C^{\infty}_c(0,+\infty)$
\begin{multline}\label{w2} D_{2\lambda^{p-1}_k}w = -
\frac{\partial}{\partial t}\left( \frac{1}{f} \frac{\partial
w}{\partial t} \right)+  \left\{ -\frac{7}{16} \frac{1}{f^3}
\left(\frac{\partial f}{\partial t}\right)^2+ \frac{1}{4}
\frac{1}{f^2} \frac{\partial^2 f}{\partial t^2} \right. \\-
\frac{1}{2} \frac{1}{f^2}\frac{\partial f}{\partial
t}\frac{(N-1+2p)}{4} \frac{1}{g}\frac{\partial g}{\partial t} +
\frac{1}{f} \frac{(N-2p+1)}{4}
\frac{(N-2p+5)}{4}\frac{1}{g^2}\left(\frac{\partial g}{\partial
t}\right)^2 \\ \left.
+\frac{1}{f}\frac{(-N+2p-1)}{4}\frac{1}{g}\frac{\partial^2
g}{\partial t^2} + \frac{\lambda^{p-1}_k}{g}\right\} w.
\end{multline}
Finally, for every $\omega \in L^2_p(M)\cap
C^{\infty}(\Lambda^p(M))$ we decompose $\omega_{2\delta}$ with
respect to an orthonormal basis of coclosed eigenforms $\left\{
\tau_{3k}\right\}_{k \in \mathbb N}$ of
$\Delta^{p-1}_{\Sphere^{N-1}}$. For every $k \in \mathbb N$ we
denote by $\lambda_k^{p-1}$ the eigenvalue corresponding to the
eigenform $\tau_{3k}$; then
$\left\{\frac{1}{\sqrt{\lambda^{p-1}_k}}d^{p-1}_{\Sphere^{N-1}}
\tau_{3k} \right\}_{k \in \mathbb N}$ is an orthonormal basis of
closed eigenforms of $\Delta^p_{\Sphere^{N-1}}$. Hence, we get the
following decomposition for any $\omega_{1d} \oplus
\omega_{2\delta}\wedge dt$
\begin{multline}\label{3}\omega_{1d}\oplus \omega_{2\delta}\wedge dt =
\oplus_k \left(\frac{1}{\sqrt{\lambda^{p-1}_k}} h_{1k}
d^{p-1}_{\Sphere^{N-1}} \tau_{3k} \oplus (-1)^p
h_{2k}\tau_{3k}\wedge dt\right),
\end{multline}
whence, by closure $${\mathcal L}_{p,3}(M)=\oplus_{k\in \mathbb
N}{\mathcal L}_{p,3,k}(M).$$ We call $p$-form of type III any
$p$-form $\omega \in L^2_p(M)$ such that $\omega \in {\mathcal
L}_{p,3,k}(M)$ for some $k\in \mathbb N$. \par A direct
computation shows that, for every $k\in \mathbb N$,
\begin{multline}\label{h1h2}
\Delta^p_{M3}\left(\frac{1}{\sqrt{\lambda_k^{p-1}}}
h_{1k}(t)d^{p-1}_{\Sphere^{N-1}} \tau_{3k} \oplus_M (-1)^p
h_{2k}(t) \tau_{3k} \wedge dt\right)\\
=\left(\Delta_{1\lambda_k^{p-1}} h_{1k} +
\frac{1}{f(t)}\frac{1}{g(t)} \frac{\partial g}{\partial
t}\sqrt{\lambda_k^{p-1}} h_{2k} \right) \left(
\frac{1}{\sqrt{\lambda_k^{p-1}}} d^{p-1}_{\Sphere^{N-1}} \tau_{3k}
\right)\\ \oplus \left(\Delta_{2\lambda_k^{p-1}} h_{2k} +
\frac{1}{g^2(t)} \frac{\partial g}{\partial t}
\sqrt{\lambda_k^{p-1}} h_{1k} \right) \left((-1)^p \tau_{3k}
\wedge dt \right);
\end{multline}
moreover, if $\omega=\frac{1}{\sqrt{\lambda}}
h_1(t)d^{p-1}_{\Sphere^{N-1}} \tau_3 \oplus_M (-1)^p h_2(t) \tau_3
\wedge dt$, then $$\|\omega\|^2_{L^2_p(M)}= \int_c^{+\infty}
g(s)^{\frac{N-2p-1}{2}}f(s)^{\frac{1}{2}} h_1(s)^2 \,ds $$ $$+
\int_c^{+\infty} g(s)^{\frac{N+1-2p}{2}}f(s)^{-\frac{1}{2}}
h_2(s)^2\,ds. $$ Hence, if we set
\begin{equation}\label{trasf12}\begin{array}{ll}
w_1(t)=&g^{\frac{N-2p-1}{4}}(t)f^{\frac{1}{4}}(t)h_1(t), \\
w_2(t)=& g^{\frac{N-2p+1}{4}}(t)f^{-\frac{1}{4}}(t) h_2(t),
\end{array}
\end{equation}
we find that $\Delta^p_{M3}$ is unitarily equivalent to the direct
sum, over $k \in \mathbb N$, of certain selfadjoint operators
$D_{3\lambda^{p-1}_k}$ in $L^2(0,+\infty)\oplus L^2(0,+\infty)$
such that $$ C^{\infty}_c(0,+\infty)\oplus
C^{\infty}_c(0,+\infty)\subseteq {\mathcal
D}(D_{3\lambda^{p-1}_k})$$ and for every $w_1 \oplus w_2 \in
C^{\infty}_c(0,+\infty)\oplus C^{\infty}_c(0,+\infty)$
\begin{multline}\label{w1w2}D_{3\lambda^{p-1}_k}(w_1\oplus w_2)=
\left(D_{1\lambda^{p-1}_k}w_1 +
g(t)^{-\frac{3}{2}}f(t)^{-\frac{1}{2}}\frac{\partial g}{\partial
t} \sqrt{\lambda^{p-1}_k}w_2 \right)  \\ \oplus
\left(D_{2\lambda^{p-1}_k} w_2 +
g(t)^{-\frac{3}{2}}f(t)^{-\frac{1}{2}}\frac{\partial g}{\partial
t} \sqrt{\lambda^{p-1}_k}w_1 \right).
\end{multline}
For $i=1,2,3$, for any $k\in \mathbb N$, we will denote by
${\mathcal T}_{p,i,k}(M)$ the unitary equivalence between
${\mathcal L}_{p,i,k}(M)$ and $L^2(0,+\infty)$ ( $
L^2(0,+\infty)\oplus L^2(0,+\infty)$ if $i=3$) given by
(\ref{trasf1}) (resp. (\ref{trasf2}), (\ref{trasf12})). \par
 We remark that even if the orthogonal decompositions depend
on the Riemannian metric (since we have to take closures in the
$L^2$-norm), the eigenvalues $\lambda_k^p$ and the eigenforms
$\tau_{1k}$ depend only on $\Sphere^{N-1}$ and hence are the same
for any choice of the functions $f(t)$ and $g(t)$. \par As a
consequence, we have:
\begin{lem}\label{XK1} Let $M$ be the $N$-dimensional
unitary ball $B(0,1)$ endowed with any complete Riemannian metric
of type
\begin{equation}ds^2= f(t)\,dt^2 +
g(t)\,d\theta_{\Sphere^{N-1}}^2,\end{equation} where $t =
\rm{settanh} (\|x\|)$, $d\theta^2_{\Sphere^{N-1}}$ is the standard
Riemannian metric on $\Sphere^{N-1}$, $f(t)>0$ and $g(t)>0$ for
any $t\in (0,+\infty)$. Moreover, let us suppose that $f$, $g$
fulfill condition (\ref{condzero}). Then for every $p\in [0,N]$,
for any $i=1,2,3$, for every $k\in \mathbb N$, the set
\begin{equation}\label{XK} X_{p,i,k}:= ({\mathcal
T}_{p,i,k}(M))(C^\infty_c(\Lambda^p(M))\cap {\mathcal
L}_{p,i,k}(M))
\end{equation}
does not depend on the choice of the functions $f$, $g$, provided
they fulfill condition (\ref{condzero}).
\end{lem}
\begin{proof}
Let $(f_1,g_1)$, $(f_2,g_2)$ be two couples of smooth positive
functions on $(0,+\infty)$ such that
\begin{enumerate}
\item the corresponding Riemannian metrics $$f_j(t)dt^2 +
g_j(t)d\theta_{\Sphere^{N-1}}, \,\,j=1,2$$ are complete on
$B(0,1)$;\\
\item for $j=1,2$, $f_j(t)=1$ and $g_j(t)=t^2$ for $t\in (0,\epsilon)$.
\end{enumerate}
For sake of simplicity, let us consider the case $i=1$ (the proofs
of the other cases are analogous). Let $\omega=h(t)\tau_{1k} \in
C^\infty_c(\Lambda^p(M))\cap {\mathcal L}_{p,1,k}(M,g_1)$; then if
we consider the differential form $\tilde \omega$ on $M$ defined
as  $$\tilde \omega
(t,\theta)=f_2(t)^{-\frac{1}{4}}g_2(t)^{-\frac{N-2p-1}{4}}f_1(t)^{\frac{1}{4}}g_1(t)^{\frac{N-2p-1}{4}}
\omega(t,\theta),$$ then $\tilde \omega \in
C^\infty_c(\Lambda^p(M)$. Moreover, it is immediate to see that
$${\mathcal T}_{p,i,k}(M,g_1)\omega={\mathcal
T}_{p,i,k}(M,g_2)\tilde \omega.$$
\end{proof}
The set $X_{p,i,k}$ defined above is a natural core for the
operator $D_{i\lambda_k}$. Namely, we have the following
characterization of ${\mathcal D}(D_{i\lambda_k})$ for any
$i=1,2,3$ and for every $k\in \mathbb N$:
\begin{lem}\label{XXK} Let $M$ be as in Lemma \ref{XK1}.
Then, for every $p\in [0,N]$, for every $i=1,2,3$, for every $k\in
\mathbb N$, the operator $D_{i\lambda_k}$ is essentially
selfadjoint on the set $X_{p,i,k}$ defined by (\ref{XK}).
\end{lem}
\begin{proof}  Since $\Delta^p_M$ is
essentially selfadjoint on $C^{\infty}_c(\Lambda^p(M))$, then, for
$i=1,2,3$, $\Delta^p_{Mi}$ is essentially selfadjoint on
$C^{\infty}_c(\Lambda^p(M))\cap {\mathcal L}_{p,i}(M)$.
Analogously, for any $i=1,2,3$ and for any $k\in \mathbb N$ the
restriction of $\Delta^p_{Mi}$ to the subspace ${\mathcal
L}_{p,i,k}(M)$ is essentially selfadjoint on
$C^{\infty}_c(\Lambda^p(M))\cap{\mathcal L}_{p,i,k}(M)$. Hence,
for every $k\in \mathbb N$ the operator $D_{i\lambda_k^p}$ is
essentially selfadjoint on the set $X_{p,i,k}$.
\end{proof}
Applying the decomposition techniques described above to the
Frie\-drichs extension $(\Delta^p_M)^F$ of the restriction of
$\Delta^p_M$ to $C^\infty_c(\Lambda^p(M\setminus B(0,c)))$ for
some arbitrarily chosen $c>0$, in \cite{Antoci2} we computed
explicitely the essential spectrum of $\Delta^p_M$ (it was shown
by Eichhorn (\cite{Eichhorn}) that $\sigma_{\rm ess}(\Delta^p_M)=
\sigma_{\rm ess}((\Delta^p_M)^F)$). Namely, we obtained the
following result:
\begin{thm}\label{spettroessenziale} Let $M$ be the unitary ball
$B(0,1)$ in $\Real^N$ endowed with a Riemannian metric $ds^2$
which, in a tubular neighbourhood of the boundary $\Sphere^{N-1}$,
is given by  $$d\sigma^2= e^{-2(a+1)t}\,dt^2 +
e^{-2bt}\,d\theta^2_{\Sphere^{N-1}},$$ where $a\leq -1$,
$t=settanh(\|x\|) $ and $d\theta^2_{\Sphere^{N-1}}$ is the
standard Riemannian metric on $\Sphere^{N-1}$. Then
\begin{enumerate}
\item if $a=-1$ and $b<0$, if $p \not= \frac{N}{2}$ $$\sigma_{\rm
ess}(\Delta^p_M)= \left[\min\left\{
\left(\frac{N-2p-1}{2}\right)^2b^2,\left(\frac{N-2p+1}{2}\right)^2b^2
\right\},+\infty\right) $$ whilst if $p=\frac{N}{2}$ $$\sigma_{\rm
ess}(\Delta^p_M)= \left\{0\right\}\cup \left[ \frac{b^2}{4},
+\infty \right) ;$$
\item if $a=-1$ and $b=0$, for every $p\in [0,N]$ $$\sigma_{\rm
ess}(\Delta^p_M)=[\overline \lambda_p,+\infty),$$ where $\overline
\lambda_p$ is the minimum between the smallest eigenvalue of
$\Delta^p_{\Sphere^{N-1}}$ and the smallest eigenvalue of
$\Delta^{p-1}_{\Sphere^{N-1}}$;
\item if $a=-1$ and $b>0$, if $1<p<N-1$ $$\sigma_{\rm ess}(\Delta^p_M)=
\emptyset,$$ whilst if $p\in \left\{0,1,N-1,N \right\}$
$$\sigma_{\rm ess}(\Delta^p_M)= \left[
\left(\frac{N-1}{2}\right)^2 b^2, +\infty\right);$$
\item if $a<-1$ and $b<0$, for every $p\in [0,N]$ $$\sigma_{\rm
ess}(\Delta^p_M)=[0,+\infty);$$
\item if $a<-1$ and $b=0$, for every $p\in [0,N]$ $$\sigma_{\rm
ess}(\Delta^p_M)=[\overline \lambda_p,+\infty),$$ where
$\overline{\lambda}$ is the minimum between the smallest
eigenvalue of $\Delta^p_{\Sphere^{N-1}}$ and the smallest
eigenvalue of $\Delta^{p-1}_{\Sphere^{N-1}}$;
\item if $a<-1$ and $b>0$, if $1<p<N-1$ $$\sigma_{\rm ess}(\Delta^p_M)=
\emptyset,$$ whilst if $p\in \left\{0,1,N-1,N\right\}$
$$\sigma_{\rm ess}(\Delta^p_M)= \left[ 0, +\infty\right).$$
\end{enumerate}
\end{thm}
As for the absolutely continuous spectrum and the singularly
continuous spectrum, in view of Lemma \ref{basic1}, for $i=1,2,3$
we have that $$\sigma_{\rm ac}(\Delta^p_{Mi})= \bigcup_{k \in
\mathbb N}\sigma_{\rm ac}(D_{i\lambda_k}) $$ and $$\sigma_{\rm
sc}(\Delta^p_{Mi})= \bigcup_{k \in \mathbb N}\sigma_{\rm
sc}(D_{i\lambda_k}) ;$$ thus, we can reduce ourselves to the
analysis of the absolutely continuous and of the singularly
continuous spectra of the selfadjoint operators
$D_{1\lambda_k^p}$, $D_{2\lambda_k^{p-1}}$ and
$D_{3\lambda_k^{p-1}}$. Since the Hodge $*$ operator maps
$p$-forms of type I isometrically onto $(N-p)$-forms of type II,
it suffices to consider the cases $i=1$ and $i=3$.
\par
Finally, let us observe that the decomposition techniques
described above work also in the case of the Euclidean space (this
will be essential in the construction of the unperturbed operators
employed in the computation of the absolutely continuous
spectrum). \par Namely,  let us consider the Euclidean space
$(\Real^N,e)$, that is, $\Real^N$ endowed with the Euclidean
metric. From now on, we will denote by $\Delta^p_e$ the
Laplace-Beltrami operator acting on $p$-forms on $(\Real^N,e)$. In
polar coordinates the Euclidean metric has the expression $$ds^2=
dr^2 + r^2 d\theta^2_{\Sphere^{N-1}},$$ where
$d\theta^2_{\Sphere^{N-1}}$ is the standard Riemannian metric on
$\Sphere^{N-1}$. Then it is possible to introduce the
decompositions $$ L^2_p(\Real^N,e)={\mathcal
L}_{p,1}(\Real^N,e)\oplus{\mathcal
L}_{p,2}(\Real^N,e)\oplus{\mathcal L}_{p,3}(\Real^N,e)$$ and, for
$i=1,2,3$, $$ {\mathcal L}_{p,i}(\Real^N,e)= \oplus_{k\in \mathbb
N}{\mathcal L}_{p,i,k}(\Real^N,e).$$ For any $k\in \mathbb N$, we
will denote by ${\mathcal T}_{p,1,k}(\Real^N,e)$ the unitary
equivalence between ${\mathcal L}_{p,1,k}(\Real^N,e)$ and $
L^2(0,+\infty)$.
\par
\section{The case $a=-1$}

Let us introduce the change of coordinates $$r:
(0,+\infty)\longrightarrow (0,+\infty),$$ $$ r(t):=
\int_0^{+\infty}\sqrt{f(s)}\,ds;$$ the Riemannian metric in the
new coordinate system $(r,\theta)$ on $(0,+\infty)\times
\Sphere^{N-1}$ is given by
\begin{equation}d\sigma^2= dr^2 + \tilde
g(r)\,d\theta^2_{\Sphere^{N-1}},\end{equation} where $$\tilde
g(r)= r^2 \quad \hbox{for}  \quad r\in (0,+\epsilon)$$ and
$$\tilde g(r)= e^{-2br}\quad \hbox{for}\quad r>\overline c=
K+\epsilon, $$ where $K=\int_{\epsilon}^c \sqrt{f(s)}\,ds$.
Applying the orthogonal decomposition of Section 4 in the new
coordinate system we find the following expressions for the
operators $D_{i\lambda_k}$: for any $w \in
C^{\infty}(0,+\infty)\cap {\mathcal D}(D_{1\lambda_k^p})$
$$D_{1\lambda_k^p}w=- \frac{\partial^2 w}{\partial r^2}+ V_1(r)w
,$$ where \begin{equation}\label{V1}V_1(r)=\begin{cases}
\left(\frac{N-2p-1}{2}\right)\left(\frac{N-2p-3}{2}\right)\frac{1}{r^2}+
\frac{\lambda_k^p}{r^2} \quad &\hbox{for $r\in (0,\epsilon)$}\\
\hbox{a smooth function}\quad &\hbox{for $r\in [\epsilon,
\overline c]$}\\ \left(\frac{N-2p-1}{2}\right)^2 b^2 + \lambda_k^p
e^{2br} \quad &\hbox{for $r>\overline c$}.
\end{cases}
\end{equation} Analogously, for any $w\in C^{\infty}(0,+\infty)\cap {\mathcal
D}(D_{2\lambda_k^{p-1}})$ $$D_{2\lambda_k^{p-1}}w=-
\frac{\partial^2 w}{\partial r^2}+ V_2(r)w ,$$ where
$$V_2(r)=\begin{cases}
\left(\frac{N-2p+1}{2}\right)\left(\frac{N-2p+3}{2}\right)\frac{1}{r^2}+
\frac{\lambda_k^{p-1}}{r^2} \quad &\hbox{for $r\in
(0,\epsilon)$}\\ \hbox{a smooth function}\quad &\hbox{for $r\in
[\epsilon, \overline c]$}\\ \left(\frac{N-2p+1}{2}\right)^2 b^2 +
\lambda_k^{p-1} e^{2br} \quad &\hbox{for $r>\overline c$}.
\end{cases} $$
Finally, for every $(w_1\oplus w_2)\in (C^\infty(0,+\infty)\oplus
C^\infty(0,+\infty))\cap {\mathcal D}(D_{3\lambda_k^{p-1}})$,
\begin{multline} D_{3\lambda^{p-1}_k}(w_1\oplus w_2)=
\left(D_{1\lambda^{p-1}_k}w_1 + V_3(r) \sqrt{\lambda^{p-1}_k}w_2
\right)  \\ \oplus \left(D_{2\lambda^{p-1}_k} w_2 + V_3(r)
\sqrt{\lambda^{p-1}_k}w_1, \right).
\end{multline}
where $$ V_3(r)=\begin{cases} \frac{2}{r^2}
 \quad &\hbox{for $r\in (0,\epsilon)$}\\
\hbox{a smooth function}\quad &\hbox{for $r\in [\epsilon,
\overline c]$}\\ -2b e^{br} \quad &\hbox{for $r>\overline c$}.
\end{cases} $$
 The behaviour of the potential at infinity depends strongly on
the sign of $b\in \Real$. Hence, we will consider separately the
cases $b<0$, $b=0$ and $b>0$.
\par \bigskip
\subsection{The case $b<0$}
We begin with the study of the absolutely continuous (and of the
singularly continuous) spectrum of the operators
$D_{1\lambda_k^p}$. To this purpose, we need some preliminary
Lemmas. The first is a classical statement
 in functional analysis:
\begin{lem}\label{confrontodom}(\cite{Reed-Simon2}) Let $A$, $C$ be symmetric operators.
Suppose that ${\mathcal
D}$ is a
linear subspace satisfying ${\mathcal D}\subseteq {\mathcal D}(A)$,
${\mathcal D}\subseteq {\mathcal D}(C)$, and that
$$\|(A-C)\varphi\|\leq a(\|A\varphi\|+\|C\varphi\|)+b\|\varphi\|$$
for all $\varphi \in {\mathcal D}$, where $0\leq a<1$, $b\geq 0$. Then
\begin{enumerate}
\item $A$ is essentially selfadjoint on ${\mathcal D}$ if and only if $C$
is essentially selfadjoint on ${\mathcal D}$;\\
\item ${\mathcal D}(\overline{A_{|{\mathcal D}}})={\mathcal
D}(\overline{C_{|{\mathcal D}}})$.
\end{enumerate}
\end{lem}
\begin{proof} For a proof see \cite{Reed-Simon2}.\end{proof}
The second Lemma is an easy generalization to the case of
differential forms of the Agmon-Kato-Kuroda Theorem (see
\cite{Reed-Simon4}).  \par We recall that a potential $V(x)$ on
$\Real^N$ is called an Agmon potential if for some $\epsilon>0$
the potential $W(x):=(1+|x|^2)^{\frac{1}{2}+\epsilon}V(x)$ is a
relatively compact perturbation of the scalar Laplacian $-\Delta$.
Moreover, it is well-known that if for some $\epsilon>0$
$(1+|x|^2)^{\frac{1}{2}+\epsilon}V(x)\in L^{\infty}(\Real^N)$ then
$V(x)$ is an Agmon potential (see \cite{Reed-Simon4}).
\begin{lem}\label{AgmonKatoKuroda} Let $V$ be an Agmon potential on $\Real^N$. If $H=\Delta^p_e+ V$,
then:
\begin{enumerate}
\item the set ${\mathcal E}_+$ of positive eigenvalues of $H$ is a
discrete subset of $(0,+\infty)$, and each eigenvalue has finite
multiplicity;\\
\item $\sigma_{\rm sc}(H)=\emptyset$;\\
\item the wave operators $W^{\pm}(H,\Delta^p_e)$ exist and are complete.
\end{enumerate}
\end{lem}
\begin{proof} For the scalar case (i.e. $p=0$) see \cite{Reed-Simon4}. For
$p>0$ the conclusion follows applying to each component the result
in the scalar case.
\end{proof}
We are now in position to prove our first result:
\begin{lem}\label{-1<01lambda} For $a=-1$, $b<0$, for $0\leq p\leq N-1$,
for every $k\in \mathbb N$
$$\sigma_{\rm ac}(D_{1\lambda_k^p})=
\left[\left(\frac{N-2p-1}{2}\right)^2 b^2,+\infty\right)\quad
\hbox{and}\quad \sigma_{\rm sc}(D_{1\lambda_k^p})= \emptyset.$$
\end{lem}
\begin{proof} We will compute the absolutely continuous
and the singularly continuous spectrum of $D_{1\lambda_k^p}$
through pertubation techniques. Since for $b<0$ we have that
$e^{2br}\rightarrow 0$ as $r\rightarrow +\infty$, it might seem
natural to apply directly the Agmon-Kato-Kuroda Theorem for
functions (see \cite{Reed-Simon4}) to the couple of operators
$(D_{1\lambda_k^p}, H)$ on the half-line $(0,+\infty)$, where
 $$H:=-\frac{\partial^2}{\partial r^2}+
\left(\frac{N-2p-1}{2}\right)^2b^2.$$ However, on one hand, the
Agmon-Kato-Kuroda Theorem holds for operators acting on the whole
$R^N$, whilst the operators $H$ and $D_{1\lambda_k^p}$ act on the
half-line. On the other hand, the potential part of the operator
$D_{1\lambda_k^p}$ has a singularity at zero. Hence, we developed
a different argument. The idea is to ``move" the problem to the
$N$-dimensional Euclidean space $(\Real^N,e)$, where the
singularity disappears.\par
 Let us consider, on $(\Real^N,e)$, the
operators $$ \tilde H_0= \Delta^p_e +
\left(\frac{N-2p-1}{2}\right)^2 b^2,$$ $$\tilde H_1= \tilde H_0 +
\tilde V(|x|), $$ where $$ \tilde V(|x|) =
\begin{cases} -\left(\frac{N-2p-1}{2}\right)^2b^2\quad \hbox{for
$|x|\in(0,\epsilon)$}\\
V_1(|x|)-\left(\left(\frac{N-2p-1}{2}\right)\left(\frac{N-2p-3}{2}\right)+
\lambda_k^p \right)\frac{1}{|x|^2}\\ -
\left(\frac{N-2p-1}{2}\right)^2b^2 \quad \hbox{for $|x|\in
[\epsilon,\overline c]$}\\
-\left(\left(\frac{N-2p-1}{2}\right)\left(\frac{N-2p-3}{2}
\right)+\lambda_k^p\right)\frac{1}{|x|^2}+\lambda_k^p
e^{2b|x|}\quad \hbox{for $|x|>\overline c$} .\end{cases} $$ Since
$\Delta^p_e$ is essentially selfadjoint on
$C^\infty_c(\Lambda^p(\Real^N,e))$, in view of Lemma
\ref{confrontodom} both $\tilde H_0$ and $\tilde H_1$ are
essentially selfadjoint on $C^{\infty}_c(\Lambda^p(\Real^N,e))$.
We denote again by $\tilde H_0$ and $\tilde H_1$ their closures.
Since an easy computation shows that for $0<\epsilon <
\frac{1}{2}$ $$(1+|x|^2)^{\frac{1}{2}+\epsilon}\tilde V(|x|)\in
L^{\infty}(\Real^N,e),$$ $\tilde V(|x|)$ is an Agmon potential on
$\Real^N$. As a consequence, Lemma \ref{AgmonKatoKuroda} implies
that
\begin{enumerate}
\item the set $\tilde{\mathcal E}$ of the eigenvalues of $\tilde H_1$
greater than $\left(\frac{N-2p-1}{2}\right)^2b^2$ is a discrete
subset of $\left(\left(\frac{N-2p-1}{2}\right)^2b^2
,+\infty\right)$, and each eigenvalue has finite multiplicity; \\
\item $\sigma_{sc}(\tilde H_1)= \emptyset$;\\
\item the wave operators $W^{\pm}(\tilde H_1, \tilde H_0)$ exist
and are complete.\end{enumerate} Now, let us consider the
restrictions $P_{|{\mathcal L}_{p,1,k}(\Real^N,e)}\tilde H_1$,
$P_{|{\mathcal L}_{p,1,k}(\Real^N,e)}\tilde H_0$ of $\tilde H_1$
and $\tilde H_0$ to ${\mathcal L}_{p,1,k}(\Real^N,e)$, and let us
apply the unitary transformation ${\mathcal
T}_{p,1,k}(\Real^N,e)$. In this way we find two operators $$H_0:=
{\mathcal T}_{p,1,k}(\Real^N,e)\circ \tilde H_0 \circ
\left({\mathcal T }_{p,1,k}(\Real^N,e)\right)^{-1}, $$ $$H_1:=
{\mathcal T}_{p,1,k}(\Real^N,e)\circ \tilde H_1 \circ
\left({\mathcal T }_{p,1,k}(\Real^N,e)\right)^{-1} ,$$ both
essentially selfadjoint on the set $X_{p,1,k}$ defined by
(\ref{XK}).\par Since a simple computation shows that for any
$w\in X_{p,1,k}$ $$H_1 w= D_{1\lambda_k^p}w,$$ in view of Lemma
\ref{XXK} we find that $H_1=D_{1\lambda^p_k}$. \par Recalling
Lemma \ref{basic1}, we find immediately that $\sigma_{\rm
sc}(D_{1\lambda_k^p})=\sigma_{\rm sc}(H_1)$ $\subseteq \sigma_{\rm
sc}(\tilde H_1)=\emptyset$. Moreover, since the projection
$P_{|{\mathcal L}_{p,1,k}(\Real^N,e)}$ commutes with both $\tilde
H_0$ and $\tilde H_1$, we find that the existence and completeness
of the wave operators $W^{\pm}(\tilde H_1,\tilde H_0)$ implies the
existence and completeness of the wave operators
$W^{\pm}(D_{1\lambda_k^p}, H_0)$. As a consequence, we have that
$$ \sigma_{\rm ac}(D_{1\lambda_k^p})= \sigma_{\rm ac}(H_0).$$
Since the spectrum of $\Delta^p_e$ is purely absolutely
continuous, equal to $[0,+\infty)$ and of constant multiplicity,
$\sigma_{\rm
ac}(H_0)=\left[\left(\frac{N-2p-1}{2}\right)^2b^2,+\infty\right)$.
This completes the proof.
\end{proof}
Hence:
\begin{prop}\label{-1<0M1} For $a=-1$, $b<0$, for $0\leq p\leq N-1$,
$$\sigma_{\rm ac}(\Delta^p_{M1})=
\left[\left(\frac{N-2p-1}{2}\right)^2
b^2,+\infty\right),\quad\hbox{and}\quad \sigma_{\rm
sc}(\Delta^p_{M1})= \emptyset.$$
\end{prop}

By duality:
\begin{prop}\label{-1<0M2}  For $a=-1$, $b<0$, for $1\leq p\leq N$,
$$\sigma_{\rm ac}(\Delta^p_{M2})=
\left[\left(\frac{N-2p+1}{2}\right)^2
b^2,+\infty\right),\quad\hbox{and}\quad \sigma_{\rm
sc}(\Delta^p_{M2})= \emptyset.$$
\end{prop}
As a consequence, since we already know from Theorem
\ref{spettroessenziale} that for $a=-1$, $b<0$, for every $p\in
[0,N]$ the essential spectrum of $\Delta^p_M$ is equal to
$\left[\min \left\{\left(\frac{N-2p-1}{2}\right)^2
b^2,\left(\frac{N-2p+1}{2}\right)^2 b^2\right\},+\infty\right),$
we can state the following:
\begin{thm}\label{-1<0}  For $a=-1$, $b<0$, for $0\leq p\leq N$,
$$\sigma_{\rm ac}(\Delta^p_M)= \left[\min
\left\{\left(\frac{N-2p-1}{2}\right)^2
b^2,\left(\frac{N-2p+1}{2}\right)^2 b^2\right\},+\infty\right),$$
$$ \sigma_{\rm sc}(\Delta^p_M)= \sigma_{\rm sc}(\Delta^p_{M3}).$$
\end{thm}

\subsection{The case $b=0$}

As in the previous case, we begin with the study of
$D_{1\lambda_k^p}$ for any $k\in \mathbb N$. If $b=0$, the
potential $V_1(r)$ in (\ref{V1}) is simply given by  $$V_1(r)=
\begin{cases}
\left(\frac{N-2p-1}{2}\right)\left(\frac{N-2p-3}{2}\right)\frac{1}{r^2}+\lambda_k^p
\frac{1}{r^2}\quad &\hbox{for $r\in (0,\epsilon)$}\\ \hbox{a
smooth function}\quad &\hbox{for $r\in [\epsilon,\overline c]$}\\
\lambda_k^p \quad &\hbox{for $r>\overline c$}
\end{cases} $$
\begin{lem}\label{-1=01lambda}  For $a=-1$, $b=0$, for $0\leq p\leq N-1$,
for every $k\in \mathbb N$
$$\sigma_{\rm ac}(D_{1\lambda_k^p})=
\left[ \lambda_k^p,+\infty\right)
\quad\hbox{and}\quad  \sigma_{\rm sc}(D_{1\lambda_k^p})= \emptyset.$$
\end{lem}
\begin{proof}
Let us consider, on $(\Real^N,e)$, the operators $$\tilde
H_0=\Delta^p_e+\lambda_k^p, $$ $$\tilde H_1=\Delta^p_e
+\lambda_k^p+ \tilde V(|x|),$$ where $$\tilde V(|x|)=\begin{cases}
-\lambda_k^p \quad &\hbox{for $|x|\in (0,\epsilon)$}\\ \hbox{a
smooth function}\quad &\hbox{for $|x|\in [\epsilon, \overline
c]$}\\ -\left(\frac{N-2p-1}{2}\right)\left(\frac{N-2p-3}{2}\right)
\frac{1}{|x|^2}-\lambda_k^p \frac{1}{|x|^2} \quad &\hbox{for
$|x|>\overline c$}.
\end{cases}$$
Again, in view of Lemma \ref{confrontodom}, both $\tilde H_1$ and
$\tilde H_0$ are essentially selfadjoint on
$C^\infty_c(\Lambda^p(\Real^N,e))$. Hence, the operators
$$H_0:={\mathcal T}_{p,1,k}(\Real^N,e)\circ \tilde H_0 \circ
\left({\mathcal T}_{p,1,k}(\Real^N,e)\right)^{-1},$$
$$H_1:={\mathcal T}_{p,1,k}(\Real^N,e)\circ \tilde H_1 \circ
\left({\mathcal T}_{p,1,k}(\Real^N,e)\right)^{-1},$$ are both
essentially selfadjoint on the set $X_{p,1,k}$. In particular, as
in the proof of Lemma \ref{-1<01lambda} we have that $H_1=
D_{1\lambda_k^p}$. Since an easy computation shows that $\tilde
V(|x|)$ is an Agmon potential on $\Real^N$ (indeed, for $0<\eps <
\frac{1}{2}$, $(1+|x|^2)^{\frac{1}{2}+\eps}\in
L^\infty(\Real^N)$), reasoning as in the proof of Lemma
\ref{-1<01lambda} we find that $\sigma_{\rm
sc}(D_{1\lambda_k^p})=\emptyset$ and $\sigma_{\rm
ac}(D_{1\lambda_k^p})=\sigma_{\rm ac}(\tilde H_0)=[\lambda_k^p,
+\infty)$.
 \end{proof} As a consequence, by Lemma \ref{basic1}, we have:
\begin{prop}\label{-1=0M1}  For $a=-1$, $b=0$, for $0\leq p\leq N-1$,
$$\sigma_{\rm ac}(\Delta^p_{M1})=
\left[\lambda_0^p,+\infty\right),$$ where $\lambda_0^p$ is the
lowest eigenvalue of $\Delta^p_{\Sphere^{N-1}}$, and $$
\sigma_{\rm sc}(\Delta^p_{M1})= \emptyset.$$
\end{prop}
By duality:
\begin{prop}\label{-1=0M2}  For $a=-1$, $b=0$, for $1\leq p\leq N$,
$$\sigma_{\rm ac}(\Delta_{M2})=
\left[\lambda_0^{p-1},+\infty\right),$$ where $\lambda_0^{p-1}$ is
the lowest eigenvalue of $\Delta^{p-1}_{\Sphere^{N-1}}$, and $$
\sigma_{\rm sc}(\Delta^p_{M2})= \emptyset.$$
\end{prop}
Since we already know from Theorem \ref{spettroessenziale} that
for $a=-1$, $b=0$ the essential spectrum of $\Delta^p_M$ is equal
to $[\overline \lambda_p, +\infty)$  for every $p\in [0,N]$, where
$\overline
\lambda_p=\min\left\{\lambda^p_0,\lambda_0^{p-1}\right\}$, we
obtain the following result:
\begin{thm}\label{-1=0}  For $a=-1$, $b=0$, for $0\leq p\leq N$,
$$\sigma_{\rm ac}(\Delta^p_M)=
\left[\overline\lambda_p,+\infty\right),$$ where
$\overline\lambda_p=
\min\left\{\lambda_0^p,\lambda_0^{p-1}\right\}$, and $$
\sigma_{\rm sc}(\Delta^p_M)= \sigma_{\rm sc}(\Delta^p_{M3}).$$
\end{thm}

\subsection{The case $b>0$}
As in the previous cases, in order to compute the absolutely
continuous spectrum of $\Delta^p_M$ it suffices to study the
absolutely continuous spectrum of $D_{1\lambda_k^p}$ for any $k\in
\mathbb N$:
\begin{lem}\label{-1>01lambda} For $a=-1$, $b>0$, for every $k\in
\mathbb N$ if $\lambda_k^p>0$
$$\sigma_{\rm ac}(D_{1\lambda_k^p})=\emptyset
\quad \hbox{and} \quad \sigma_{\rm sc}(D_{1\lambda_k^p})=\emptyset,$$
whilst if
$\lambda_k^p=0$ $$\sigma_{\rm
ac}(D_{1\lambda_k^p})=\left[\left(\frac{N-1}{2}\right)^2b^2,+\infty\right)
\quad \hbox{and} \quad \sigma_{\rm
sc}(D_{1\lambda_k^p})=\emptyset.$$
\end{lem}
\begin{proof} It was proved in \cite{Antoci2} that for $a=-1$, $b>0$, if
$\lambda_k^p>0$ then $\sigma_{\rm
ess}(D_{1\lambda_k^p})=\emptyset$; as a consequence, in this case
$\sigma_{\rm ac}(D_{1\lambda_k^p})=\sigma_{\rm
sc}(D_{1\lambda_k^p})=\emptyset$.\par If, on the contrary,
$\lambda_k^p=0$, we have that $V_1(r)$ is simply  $$ V_1(r)=
\begin{cases}
\left(\frac{N-2p-1}{2}\right)\left(\frac{N-2p-3}{2}\right)\frac{1}{r^2}
\quad &\hbox{for $r\in (0,\epsilon)$}\\ \hbox{a smooth
function}\quad &\hbox{for $r\in [\epsilon, \overline c]$}\\
\left(\frac{N-2p-1}{2}\right)^2 b^2\quad &\hbox{for $r> \overline
c$}.
\end{cases}$$
Let us consider, on $(\Real^N,e)$, the operators $$\tilde
H_0=\Delta^p_e+\left(\frac{N-2p-1}{2}\right)^2b^2, $$ $$\tilde
H_1=\Delta^p_e +\left(\frac{N-2p-1}{2}\right)^2b^2+ \tilde
V(|x|),$$ where $$\tilde V(|x|)=\begin{cases} -
\left(\frac{N-2p-1}{2}\right)^2b^2 \quad &\hbox{for $|x|\in
(0,\epsilon)$}\\ \hbox{a smooth function}\quad &\hbox{for $|x|\in
[\epsilon, \overline c]$}\\
-\left(\frac{N-2p-1}{2}\right)\left(\frac{N-2p-3}{2}\right)
\frac{1}{|x|^2} \quad &\hbox{for $|x|>\overline c$}.
\end{cases}$$
Again, in view of Lemma \ref{confrontodom}, both $\tilde H_1$ and
$\tilde H_0$ are essentially selfadjoint on
$C^\infty_c(\Lambda^p(\Real^N,e))$. Hence, the operators
$$H_0:={\mathcal T}_{p,1,k}(\Real^N,e)\circ \tilde H_0 \circ
\left({\mathcal T}_{p,1,k}(\Real^N,e)\right)^{-1},$$
$$H_1:={\mathcal T}_{p,1,k}(\Real^N,e)\circ \tilde H_1 \circ
\left({\mathcal T}_{p,1,k}(\Real^N,e)\right)^{-1},$$ are both
essentially selfadjoint on the set $X_{p,1,k}$. This fact, jointly
with a simple computation, shows that $H_1= D_{1\lambda_k^p}$.
\par Since $\tilde V(|x|)$ is an Agmon potential on $\Real^N$,
reasoning as in the proof of Lemma \ref{-1<01lambda} we find that
$\sigma_{\rm sc}(D_{1\lambda_k^p})=\emptyset$ and $\sigma_{\rm
ac}(D_{1\lambda_k^p})=\sigma_{\rm ac}(\tilde
H_0)=[\left(\frac{N-2p-1}{2}\right)^2b^2, +\infty)$.
\end{proof}
Now, it is well-known that on $\Sphere^{N-1}$ we can have $\lambda_k^p=0$
(that is, there exist harmonic $p$-forms) only for $p=0$ or for $p=N-1$.
\par Hence:
\begin{prop}\label{-1>0M1} For $a=-1$, $b>0$, if $0<p<N-1$ $$\sigma_{\rm
ac}(\Delta^p_{M1})=\emptyset \quad \hbox{and} \quad \sigma_{\rm
sc}(\Delta^p_{M1})=\emptyset , $$ whilst if $p\in \left\{0,N-1
\right\}$ $$\sigma_{\rm
ac}(\Delta^p_{M1})=\left[\left(\frac{N-1}{2}\right)^2b^2,+\infty\right)
\quad \hbox{and} \quad \sigma_{\rm sc}(\Delta^p_{M1})=
\emptyset.$$
\end{prop}
By duality:
\begin{prop}\label{-1>0M2} For $a=-1$, $b>0$, if $1<p<N$
$$\sigma_{\rm ac}(\Delta^p_{M2})=\emptyset \quad \hbox{and} \quad
\sigma_{\rm sc}(\Delta^p_{M2})=\emptyset , $$ whilst if $p\in
\left\{1,N\right\}$ $$\sigma_{\rm
ac}(\Delta^p_{M1})=\left[\left(\frac{N-1}{2}\right)^2b^2,+\infty\right)
\quad\hbox{and} \quad \sigma_{\rm sc}(\Delta^p_{M1})= \emptyset.$$
\end{prop}
As a consequence, in view of Theorem \ref{spettroessenziale}, we
have the following result:
\begin{thm}\label{-1>0} For $a=-1$, $b>0$, if $1<p<N-1$
$$ \sigma_{\rm ac}(\Delta^p_M)= \emptyset \quad \hbox{and} \quad
\sigma_{\rm sc}(\Delta^p_M)=\sigma_{\rm sc}(\Delta^p_{M3}),$$
whilst if $p \in \left\{ 0,1,N-1,N \right\}$ $$\sigma_{\rm
ac}(\Delta^p_M)=\left[\left(\frac{N-1}{2}\right)^2b^2,+\infty\right)
\hbox{and}\quad \sigma_{\rm sc}(\Delta^p_M)=\sigma_{\rm
sc}(\Delta^p_{M3}).$$
\end{thm}

\section{The case $a<-1$}

As in the previous section, we introduce the change of coordinates
$$r: (0,+\infty)\longrightarrow (0,+\infty),$$ $$ r(t):=
\int_0^{+\infty}\sqrt{f(s)}\,ds;$$ the Riemannian metric in the
new coordinate system $(r,\theta)$ on $(0,+\infty)\times
\Sphere^{N-1}$ is given by
\begin{equation}d\sigma^2= dr^2 + \tilde
g(r)\,d\theta^2_{\Sphere^{N-1}},\end{equation} where $$\tilde
g(r)= r^2 \quad \hbox{for}  \quad r\in (0,+\epsilon)$$ and
$$\tilde g(r)=
|a+1|^{-\frac{2b}{|a+1|}}(r-c_1)^{-\frac{2b}{|a+1|}}\quad
\hbox{for}\quad r>\overline c=K+\epsilon, $$ where
$K=\int_{\epsilon}^c \sqrt{f(s)}\,ds$ and $c_1=K+\epsilon -
\frac{e^{|a+1|c}}{|a+1|}>0$. Applying the orthogonal decomposition
of Section 4 in the new coordinate system we find the following
expression for the operators $D_{i\lambda_k}$: for any $w \in
C^{\infty}(0,+\infty)\cap{\mathcal D}(D_{1\lambda_k^p})$
$$D_{1\lambda_k^p}w= -\frac{\partial^2 w}{\partial r^2}+ V_1(r)w
,$$ where $$V_1(r)=\begin{cases}
\left(\frac{N-2p-1}{2}\right)\left(\frac{N-2p-3}{2}\right)\frac{1}{r^2}+
\frac{\lambda_k^p}{r^2} \quad &\hbox{for $r\in (0,\epsilon)$}\\
\hbox{a smooth function}\quad &\hbox{for $r\in [\epsilon,
\overline c]$}\\ \tilde K_1 (r-c_1)^{-2}+ \lambda_k^p
|a+1|^{\frac{2b}{|a+1|}}(r-c_1)^{\frac{2b}{|a+1|}} \quad
&\hbox{for $r>\overline c$},
\end{cases}
$$ where $$\tilde K_1=\left(\frac{N-2p-1}{2}\right)^2
\frac{b^2}{|a+1|^2}+ \frac{N-2p-1}{2}\frac{b}{|a+1|}.$$
Analogously, for any $w \in C^{\infty}(0,+\infty)\cap{\mathcal
D}(D_{2\lambda_k^{p-1}})$ $$D_{2\lambda_k^{p-1}}w=
-\frac{\partial^2 w}{\partial r^2}+ V_2(r)w ,$$ where
$$V_2(r)=\begin{cases}
\left(\frac{N-2p+1}{2}\right)\left(\frac{N-2p+3}{2}\right)\frac{1}{r^2}+
\frac{\lambda_k^{p-1}}{r^2} \quad &\hbox{for $r\in
(0,\epsilon)$}\\ \hbox{a smooth function}\quad &\hbox{for $r\in
[\epsilon, \overline c]$}\\ \tilde K_2 (r-c_1)^{-2}+
\lambda_k^{p-1} |a+1|^{\frac{2b}{|a+1|}}(r-c_1)^{\frac{2b}{|a+1|}}
\quad &\hbox{for $r>\overline c$},
\end{cases}
$$ where $$\tilde K_2=\left(\frac{N-2p+1}{2}\right)^2
\frac{b^2}{|a+1|^2}+ \frac{N-2p+1}{2}\frac{b}{|a+1|}.$$ Finally,
for every $(w_1\oplus w_2)\in (C^\infty(0,+\infty)\oplus
C^\infty(0,+\infty))\cap {\mathcal D}(D_{3\lambda_k^{p-1}})$,
\begin{multline} D_{3\lambda^{p-1}_k}(w_1\oplus w_2)=
\left(D_{1\lambda^{p-1}_k}w_1 + V_3(r) \sqrt{\lambda^{p-1}_k}w_2
\right)  \\ \oplus \left(D_{2\lambda^{p-1}_k} w_2 + V_3(r)
\sqrt{\lambda^{p-1}_k}w_1, \right).
\end{multline}
where $$ V_3(r)=\begin{cases} \frac{2}{r^2}
 \quad &\hbox{for $r\in (0,\epsilon)$}\\
\hbox{a smooth function}\quad &\hbox{for $r\in [\epsilon,
\overline c]$}\\|a+1|^{\frac{b}{|a+1|}}(r-c_1)^{\frac{b}{|a+1|}-1}
\quad &\hbox{for $r>\overline c$}.
\end{cases} $$
As in the previous section, the behaviour of the potential at
$+\infty$ depends strongly on the sign of $b\in \Real$, thus we
will consider separately the cases $b<0$, $b=0$ and $b>0$.
\subsection{The case $b<0$}
Let us begin with the study of the spectrum of $D_{1\lambda_k^p}$
for any $k\in \mathbb N$ . To this purpose, let us introduce the
following Theorem, which is an easy generalization to the case of
$p$-forms of a result due to Lavine (see \cite{Lavine}):
\begin{thm}\label{Lavine} Let $\tilde V$ be a multiplication operator
acting on $L^2_p(\Real^N,e)$, where $$ \tilde V (x)=
V_\alpha(x)+V_\beta(x),$$ with
\begin{enumerate}
\item $V_\alpha \in C^1(\Real^N)$,
\item $\lim_{|x|\rightarrow +\infty}V_\alpha(x)=0$
\item $|\frac{\partial V_\alpha}{\partial r}|\leq
c(1+r)^{-\gamma}$ for some $\gamma >1$ (here $r=|x|$),
\item $V_\beta(x)= (1+|x|)^{-\gamma}(f_p+f_{\infty})$ for some $\gamma >1$,
$f_{\infty}\in L^{\infty}(\Real^N)$, $f_p\in L^p(\Real^N)$ for
$p>\max(\frac{N}{2}, 1)$. \end{enumerate} Then there exists a
unique selfadjoint operator $H$ with ${\mathcal D}(H)$ $\subseteq
\\{\mathcal D}((\Delta^p_e)^{\frac{1}{2}})$ such that for every
$\omega \in {\mathcal D}(H)$ $$\langle H \omega,
\omega\rangle_{L^2_p(\Real^N,e)}= \sum_{i,j=1}^N \int_{\Real^N}
\left(\frac{\partial \omega_i}{\partial x_j}\right)^2 \,dx+
\int_{\Real^N} \tilde V(x) |\omega(x)|^2\,dx .$$ The positive
eigenvalues of $H$ have finite multiplicity and can accumulate
only at $0$. Moreover, $${\mathcal H}_{\rm ac}(H)=\left({\mathcal
H}_p(H)\right)^{\perp}.$$
\end{thm}
\begin{proof} For the scalar case (i.e. $p=0$) see \cite{Lavine}. For
$p>0$ the assert follows applying to each component the result in the
scalar case. \end{proof}
\begin{rem} Under the assumptions of Theorem \ref{Lavine}, we do not get
the existence and completeness of the wave operators
$W^{\pm}(H,\Delta^p_e)$ (indeed, for certain potentials they might
not exist, as shown in \cite{Dollard}).\end{rem} We are now in
position to prove the following
\begin{lem}\label{<-1<01lambda} For $a=-1$, $b<0$, for $0\leq p\leq N-1$,
for every $k\in \mathbb N$
$$\sigma_{\rm ac}(D_{1\lambda_k^p})=
\left[0,+\infty\right)\quad\hbox{and}\quad
\sigma_{\rm sc}(D_{1\lambda_k^p})= \emptyset.$$
\end{lem}
\begin{proof}
Let us consider, on the Euclidean space $(\Real^N,e)$, the
operators $$\tilde H_0:=\Delta^p_e,$$ $$\tilde H_1:= \Delta^p_e+
\tilde V(|x|),$$ where $$\tilde V(|x| )=\begin{cases} 0 \quad
\hbox{for $ |x|\in (0,\epsilon)$}\\ \hbox{a smooth function}\quad
\hbox{for $|x| \in [\epsilon, \overline c]$}\\ \tilde K_1 (|
x|-c_1)^{-2}+\lambda_k^p|a+1|^{-\frac{2|b|}{|a+1|}}(r-c_1)^{-\frac{2|b|}
{|a+1|}}\\ -
\left(\frac{N-2p-1}{2}\right)\left(\frac{N-2p-3}{2}\right)\frac{1}{|
x|^2} - \lambda_k^p \frac{1}{| x|^2}\quad \hbox{for $|x|>\overline
c $}.\end{cases}$$ Since $\tilde V(|x|)$ is bounded, we have that
$\tilde H_1$ is essentially selfadjoint on
$C^\infty_c(\Lambda^p(\Real^N,e))$. Hence, the operators $$
H_0:={\mathcal T}_{p,1,k}(\Real^N,e)\circ \tilde H_0 \circ
({\mathcal T}_{p,1,k}(\Real^N,e))^{-1},$$ $$ H_1:={\mathcal
T}_{p,1,k}(\Real^N,e)\circ \tilde H_1 \circ ({\mathcal
T}_{p,1,k}(\Real^N,e))^{-1},$$ are both essentially selfadjoint on
the set $X_{p,1,k}$. Since $D_{1\lambda_k^p}$ is essentially
selfadjoint on $X_{p,1,k}$ and $D_{1\lambda_k^p}w=H_1 w$ for every
$w\in X_{p,1,k}$, we have that $H_1=D_{1\lambda_k^p}$.
\par
Now, $\tilde V(|x|)$ is not an Agmon potential for any possible
value of $a<-1$, $b<0$. \par  If $|b|>\frac{|a+1|}{2}$, then for
$0< \epsilon < \min\left\{\frac{1}{2}, \frac{1}{2}
\left(\frac{2|b|}{|a+1|}-1\right)\right\}$ we have that
$$(1+|x|^2)^{\frac{1}{2}+\epsilon}\tilde V(|x|) \in
L^{\infty}(\Real^N), $$ hence $\tilde V(|x|)$ is an Agmon
potential on $\Real^N$. As a consequence, following the argument
of Lemma \ref{-1<01lambda} we find that for $|b|>\frac{|a+1|}{2}$
$\sigma_{\rm ac}(D_{1\lambda_k^p})=[0,+\infty)$ and $\sigma_{\rm
sc}(D_{1\lambda_k^p})=\emptyset$.\par If, on the contrary,
$|b|\leq \frac{|a+1|}{2}$, $\tilde V(|x|)$ is no more an Agmon
potential; however, $\tilde V(|x|)$ fulfills the assumptions of
Theorem \ref{Lavine}. Indeed, $\tilde V(|x|)$ can be written as
$$\tilde V(|x|)=V_\alpha(|x|)+V_\beta(|x|),$$ where
$$V_\alpha(|x|)=V_\beta(|x|)=0 $$ for $|x|\in (0,\epsilon)$,
whilst for $|x|>\overline c$ $$V_\alpha(|x|)= \lambda_k^p
|a+1|^{-\frac{2|b|}{|a+1|}}(|x|-c_1)^{-\frac{2|b|}{|a+1|}}$$ and
$$V_\beta(|x|)=\tilde K_1
(|x|-c_1)^{-2}-\left(\frac{N-2p-1}{2}\frac{N-2p-3}{2}
+\lambda_k^p\right)\frac{1}{|x|^2}.$$ It is immediate to see that
$V_\alpha\in C^1(\Real^N)$, $V_\alpha(|x|)\rightarrow 0$ as
$|x|\rightarrow +\infty$ and $$|\frac{\partial V_\alpha}{\partial
r}|\leq C (1+r)^{-\left(\frac{2|b|}{|a+1|}+1\right)}$$ for some
positive constant $C$.\par Moreover, for $\eps <1$
$$(1+|x|)^{1+\eps}V_2(|x|)\in L^{\infty}(\Real^N).$$ As a
consequence, by Theorem \ref{Lavine} $${\mathcal H}_{\rm
ac}(\tilde H_1)=\left({\mathcal H}_p(\tilde H_1)\right)^\perp ;$$
moreover, the positive eigenvalues of $\tilde H_1$ have finite
multiplicity and can accumulate only at $0$. These facts hold also
for the restriction of $\tilde H_1$ to the subspace ${\mathcal
L}_{p,1,k}(\Real^N,e)$. Hence, we find that, for every $k\in
\mathbb N$,
\begin{equation}\label{perp} {\mathcal H}_{\rm
ac}(D_{1\lambda_k^p})=\left({\mathcal
H}_p(D_{1\lambda_k^p})\right)^\perp;\end{equation} moreover, for
every $k\in \mathbb N$ the positive eigenvalues of
$D_{1\lambda_k^p}$ have finite multiplicity and can accumulate
only at $0$. \par From (\ref{perp}) we immediately get ${\mathcal
H}_{\rm sc}(D_{1\lambda_k^p})=\emptyset$, whence $$\sigma_{\rm
sc}(D_{1\lambda_k^p})=\emptyset.$$ As for the absolutely
continuous spectrum, since $$\sigma(D_{1\lambda_k^p})=\sigma_{\rm
ac}(D_{1\lambda_k^p})\cup \sigma_p(D_{1\lambda_k^p})$$ and, in
view of Theorem 5.1 in \cite{Antoci2}, $$ \sigma
(D_{1\lambda_k^p})=[0,+\infty),$$ we find $$[0,+\infty)\setminus
\sigma_p(D_{1\lambda_k^p})\subseteq \sigma_{\rm
ac}(D_{1\lambda_k^p}),$$ whence $$\sigma_{\rm
ac}(D_{1\lambda_k^p})=[0,+\infty).$$
\end{proof}
Hence, by Lemma \ref{basic1},
\begin{prop}\label{<-1<0M1} For $a<-1$, $b<0$, for $0\leq p\leq N-1$,
$$\sigma_{\rm ac}(\Delta^p_{M1})= \left[0,+\infty\right)\quad
\hbox{and}\quad \sigma_{\rm sc}(\Delta^p_{M1})= \emptyset.$$
\end{prop}
By duality:
\begin{prop}\label{<-1<0M2} For $a<-1$, $b<0$, for $1\leq p\leq N$,
$$\sigma_{\rm ac}(\Delta^p_{M2})= \left[0,+\infty\right)\quad
\hbox{and}\quad \sigma_{\rm sc}(\Delta^p_{M2})= \emptyset.$$
\end{prop}
As a consequence, since from Theorem \ref{spettroessenziale} we
already know that for $a<-1$, $b<0$ the essential spectrum of
$\Delta^p_M$ is equal to $[0,+\infty)$ for every $p\in [0,N]$, we
can state the following
\begin{thm}\label{<-1<0} For $a<-1$, $b<0$, for $0\leq p\leq N$,
$$\sigma_{\rm ac}(\Delta^p_M)=
\left[0,+\infty\right)\quad\hbox{and}\quad \sigma_{\rm
sc}(\Delta^p_M)= \sigma_{\rm sc}(\Delta^p_{M3}).$$
\end{thm}

\subsection{The case $b=0$}
First of all, we study the spectral properties of
$D_{1\lambda_k^p}$ for every $k\in \mathbb N$.
\begin{lem}\label{<-1=01lambda} For $a<-1$, $b=0$, for $0\leq p\leq N-1$,
for every $k\in \mathbb N$
$$\sigma_{\rm ac}(D_{1\lambda_k^p})=
\left[ \lambda_k^p,+\infty\right)\quad\hbox{and}\quad  \sigma_{\rm
sc}(D_{1\lambda_k^p})= \emptyset.$$
\end{lem}
\begin{proof} For $b=0$, the potential $V_1(r)$ is simply given by
$$ V_1(r)=
\begin{cases}\left(\frac{N-2p-1}{2}\right)\left(\frac{N-2p-3}{2}\right)
\frac{1}{r^2}+ \lambda_k^p \frac{1}{r^2} \quad &\hbox{for $r\in
(0,\epsilon)$}\\ \hbox{a smooth function} \quad &\hbox{for $r\in
[\epsilon, \overline c]$}\\ \tilde K_1 (r-c_1)^{-2}+\lambda_k^p
\quad &\hbox{for $r>\overline c$} .\end{cases}$$ Let us consider,
on the Euclidean space $(\Real^N,e)$, the operators $$\tilde H_0
:= \Delta^p_e +\lambda_k^p,$$ $$\tilde H_1:=\Delta^p_e
+\lambda_k^p + \tilde V(|x|),$$ where $$\tilde
V(|x|)=\begin{cases} - \lambda_k^p \quad \hbox{for $|x|\in
(0,\epsilon)$}\\ \hbox{a smooth function}\quad \hbox{for $|x|\in
[\epsilon,\overline c]$}\\ \tilde K_1 (|x|-c_1)^{-2}
-\left(\frac{N-2p-3}{2}\frac{N-2p-1}{2}+ \lambda_k^p\right)
\frac{1}{|x|^2} \quad \hbox{for $|x|>\overline c$} .\end{cases}$$
Since $\tilde V(|x|)$ is an Agmon potential on $\Real^N$,
following the argument of Lemma \ref{-1<01lambda} we find that
$\sigma_{\rm sc}(D_{1\lambda_k^p})=\emptyset$ and $\sigma_{\rm
ac}(D_{1\lambda_k^p})= [\lambda_k^p,+\infty)$ for every $k\in
\mathbb N$. \end{proof} Hence, by Lemma \ref{basic1}:
\begin{prop}\label{<-1=0M1} For $a<-1$, $b=0$, for $0\leq p\leq N-1$,
$$\sigma_{\rm ac}(\Delta^p_{M1})=
\left[\lambda_0^p,+\infty\right),$$ where $\lambda_0^p$ is the
lowest eigenvalue of $\Delta^p_{\Sphere^{N-1}}$ on $p$-forms, and
$$ \sigma_{\rm sc}(\Delta^p_{M1})= \emptyset.$$
\end{prop}
By duality:
\begin{prop}\label{<-1=0M2}
For $a<-1$, $b=0$, for $1\leq p\leq N$, $$\sigma_{\rm
ac}(\Delta^p_{M2})= \left[\lambda_0^{p-1},+\infty\right),$$ where
$\lambda_0^{p-1}$ is the lowest eigenvalue of
$\Delta^p_{\Sphere^{N-1}}$ on $(p-1)$-forms, and $$ \sigma_{\rm
sc}(\Delta^p_{M2})= \emptyset.$$
\end{prop}
As a consequence, since we know from Theorem
\ref{spettroessenziale} that for $a<-1$, $b=0$, for every $p\in
[0,N]$ the essential spectrum of $\Delta^p_M$ is equal to
$[\overline \lambda_p,+\infty)$, where $\overline \lambda_p= \min
\left\{ \lambda_0^p, \lambda_0^{p-1} \right\}$, we find the
following result:
\begin{thm}\label{<-1=0}  For $a<-1$, $b=0$, for $0\leq p\leq N$,
$$\sigma_{\rm ac}(\Delta^p_M)=
\left[\overline\lambda,+\infty\right),$$ where $\overline\lambda=
\min\left\{\lambda_0^p,\lambda_0^{p-1}\right\}$, and $$
\sigma_{\rm sc}(\Delta^p_M)= \sigma_{\rm sc}(\Delta^p_{M3}).$$
\end{thm}

\subsection{The case $b>0$}
As in the previous cases, in order to compute the absolutely
continuous spectrum of $\Delta^p_M$ it suffices to study the
absolutely continuous spectrum of $D_{1\lambda^p_k}$ for every
$k\in \mathbb N$:
\begin{lem}\label{<-1>01lambda}  For $a<-1$, $b>0$, for every $k\in
\mathbb N$ if $\lambda_k^p>0$
$$\sigma_{\rm ac}(D_{1\lambda_k^p})=\emptyset
\quad \hbox{and} \quad \sigma_{\rm sc}(D_{1\lambda_k^p})=\emptyset,$$
whilst if
$\lambda_k^p=0$
$$\sigma_{\rm ac}(D_{1\lambda_k^p})=[0,+\infty)\quad \hbox{and} \quad
\sigma_{\rm sc}(D_{1\lambda_k^p})=\emptyset.$$
\end{lem}
\begin{proof}
It was proved in \cite{Antoci2} that for $a<-1$, $b>0$, if
$\lambda_k^p>0$ then $\sigma_{\rm
ess}(D_{1\lambda_k^p})=\emptyset$; thus, if $\lambda_k^p>0$, then
$\sigma_{\rm ac}(D_{1\lambda_k^p})=\sigma_{\rm
sc}(D_{1\lambda_k^p})=\emptyset$.\par If, on the contrary,
$\lambda_k^p=0$, then $V_1(r)$ is simply given by  $$ V_1(r)=
\begin{cases}
\left(\frac{N-2p-1}{2}\right)\left(\frac{N-2p-3}{2}\right)\frac{1}{r^2}
\quad &\hbox{for $r\in (0,\epsilon)$}\\ \hbox{a smooth
function}\quad &\hbox{for $r\in [\epsilon, \overline c]$}\\ \tilde
K_1 (r-c_1)^{-2}\quad &\hbox{for $r> \overline c$}.
\end{cases}$$ Let us consider, on $(\Real^N,e)$, the operators
$$\tilde H_0:= \Delta^p_e,$$ $$\tilde H_1:= \Delta^p_e + \tilde
V(|x|),$$ where $$ \tilde V(|x|)=\begin{cases} 0 \quad \hbox{for
$|x|\in (0,\epsilon)$}\\ \hbox{a smooth function}\quad \hbox{for
$|x|\in [\epsilon, \overline c]$}\\ -
\left(\frac{N-2p-1}{2}\right)\left(\frac{N-2p-3}{2}\right)\frac{1}{|x|^2}+
\tilde K_1 (|x|-c_1)^{-2} \quad \hbox{for $|x|>\overline c
$}.\end{cases}$$ Since $\tilde V(|x|)$ is an Agmon potential on
$\Real^N$, following the argument of Lemma \ref{-1<01lambda} we
find that if $\lambda_k=0$ then $\sigma_{\rm
sc}(D_{1\lambda_k^p})=\emptyset$ and $\sigma_{\rm
ac}(D_{1\lambda_k^p})= \left[0, +\infty\right).$ This completes
the proof.
\end{proof}
Hence, by Lemma \ref{basic1}:
\begin{prop}\label{<-1>0M1}  For $a<-1$, $b>0$, if $0<p<N-1$ $$\sigma_{\rm
ac}(\Delta^p_{M1})=\emptyset \quad \hbox{and} \quad \sigma_{\rm
sc}(\Delta^p_{M1})=\emptyset , $$ whilst if $p\in \left\{0,N-1
\right\}$ $$\sigma_{\rm ac}(\Delta^p_{M1})=[0,+\infty) \quad
\hbox{and} \quad \sigma_{\rm sc}(\Delta^p_{M1})= \emptyset.$$
\end{prop}
By duality:
\begin{prop}\label{<-1>0M2} For $a<-1$, $b>0$, if $1<p<N$ $$\sigma_{\rm
ac}(\Delta^p_{M2})=\emptyset \quad \hbox{and} \quad \sigma_{\rm
sc}(\Delta^p_{M2})=\emptyset , $$ whilst if $p\in \left\{1,N
\right\}$ $$\sigma_{\rm ac}(\Delta^p_{M2})=[0,+\infty) \quad
\hbox{and} \quad \sigma_{\rm sc}(\Delta^p_{M2})= \emptyset.$$
\end{prop}
As a consequence, in view of Theorem \ref{spettroessenziale}, we
get the following result:
\begin{thm}\label{<-1>0}  For $a<-1$, $b>0$, if $1<p<N-1$
$$ \sigma_{\rm ac}(\Delta^p_M)= \emptyset \quad \hbox{and} \quad
\sigma_{\rm sc}(\Delta^p_M)=\sigma_{\rm sc}(\Delta^p_{M3}),$$
whilst if $p \in \left\{ 0,1,N-1,N \right\}$ $$\sigma_{\rm
ac}(\Delta^p_M)=[0,+\infty)\quad \hbox{and} \quad \sigma_{\rm
sc}(\Delta^p_M)=\sigma_{\rm sc}(\Delta^p_{M3}).$$
\end{thm}

%----------------------------------------------------------------

\end{document}